\documentclass{siamart250211}
\newsiamremark{example}{Example}

\usepackage{import}
\usepackage{xcolor}
\definecolor{uqp}{RGB}{152,24,147}
\definecolor{oxb}{RGB}{0,33,71}
\definecolor{BurntOrange}{RGB}{203,96,21}

\usepackage{mathtools}
\usepackage{amsfonts}
\usepackage{bm}  
\usepackage{mleftright} \mleftright{}  
\usepackage{braket}  
\usepackage{siunitx}  

\let\e\varepsilon{}               

\let\N\Natural{}
\let\Q\Rational{}
\let\R\Real{}
\let\Z\Integer{}

\newcommand{\ball}{B}

\newcommand{\Lip}{\mathrm{Lip}}

\newcommand{\deq}{\coloneqq}	

\DeclarePairedDelimiter{\abs}{\lvert}{\rvert}
\DeclarePairedDelimiter{\norm}{\lVert}{\rVert}
\DeclarePairedDelimiter{\paren}{\lparen}{\rparen}

\DeclarePairedDelimiterX{\innerprod}[2]{\langle}{\rangle}{#1,#2}

\newcommand{\Abs}[1]{\abs*{#1}}

\newcommand{\Paren}[1]{\paren*{#1}}

\newcommand{\Innerprod}[2]{\innerprod*{#1}{#2}}

\newcommand{\vek}[1]{\bm{#1}}
\newcommand{\mat}[1]{\bm{#1}}

\usepackage{booktabs}  
\usepackage{multirow}
\usepackage{array} 
\usepackage{tabularx} 
\usepackage{rotating} 
\usepackage{makecell} 
\usepackage{siunitx} 

\usepackage[shortlabels]{enumitem}

\usepackage{pgf}

\usepackage[linesnumbered, lined, ruled, algosection, algo2e]{algorithm2e}
\DontPrintSemicolon
\SetAlgoHangIndent{0pt}
\SetKwInput{KwParameters}{Parameters}
\SetKwInput{KwDefaults}{Defaults}
\SetKwInput{KwNote}{Note}
\crefname{algocf}{Algorithm}{Algorithms}

\ifpdf{}
\hypersetup{
  pdftitle={Local Convergence of Adaptively Regularized Tensor Methods},
  pdfauthor={Karl Welzel, Yang Liu, Raphael A. Hauser, and Coralia Cartis}
}
\fi

\headers{Local Convergence of Adaptively Regularized Tensor Methods}{Welzel, Liu, Hauser, and Cartis}

\title{Local Convergence of Adaptively Regularized Tensor Methods\thanks{Received by the editors DATE.\@
\funding{This work is supported by the Hong Kong Innovation and Technology Commission (InnoHK Project CIMDA)}}}

\author{Karl Welzel\thanks{Mathematical Institute, University of Oxford, Woodstock Road, Oxford OX2 6GG, United Kingdom
  (\email{karl.welzel@maths.ox.ac.uk}, \email{yang.liu@maths.ox.ac.uk}, \email{raphael.hauser@maths.ox.ac.uk}, \email{coralia.cartis@maths.ox.ac.uk}).}
\and Yang Liu\footnotemark[2] \and Raphael A. Hauser\footnotemark[2] \and Coralia Cartis\footnotemark[2]}

\begin{document}

\maketitle

\begin{abstract}
    Optimization methods that make use of derivatives of the objective function up to order \(p > 2\) are called tensor methods.
    Among them, ones that minimize a regularized \(p\)th-order Taylor expansion at each step have been shown to possess optimal global complexity which improves as \(p\) increases.
    The local convergence of such optimization algorithms on functions that have Lipschitz continuous \(p\)th derivatives and are uniformly convex of order \(q\) has been studied by Doikov and Nesterov [\textit{Math. Program.}, 193 (2022), pp. 315--336].
    We extend these local convergence results to locally uniformly convex functions and fully adaptive methods, which do not need knowledge of the Lipschitz constant, thus providing the first sharp local rates for AR\(p\).
    We discuss the surprising new challenges encountered by nonconvex local models and non-unique model minimizers.
    For \(p > 2\), our examples show that in particular when using the global minimizer of the subproblem, even asymptotically not all iterations need to be successful.
    Only if the ``right'' local model minimizer is used, the \(p/(q-1)\)th-order local convergence from the non-adaptive case is preserved for \(p > q-1\), otherwise the superlinear rate can degrade.
    We thus confirm that adaptive higher-order methods achieve superlinear convergence for certain degenerate problems as long as \(p\) is large enough and provide sharp bounds on the order of convergence one can expect in the limit.
\end{abstract}



\section{Introduction}

Tensor methods are a class of methods for unconstrained continuous optimization that can incorporate derivatives of order \(1\) to \(p\) of the objective function \(f \in C^p(\R^n)\) in order to find a local minimizer of \(f\).
The most important method in this regard is the AR\(p\) method.
It works by constructing the local model at each iterate $\vek{x}_k$ as the sum of the \(p\)th-order Taylor expansion at \(\vek{x}_k\) and a regularization term \(\sigma_k \norm{\vek{y} - \vek{x}_k}^{p+1}\).
This model is then minimized to find a potential new iterate \(\vek{y}_k\).
If the function decrease at the new iterate is sufficient, the iteration is called ``successful'', the iterate is accepted (\(\vek{x}_{k+1} = \vek{y}_k\)) and the regularization parameter \(\sigma_k\) is decreased; otherwise the iteration is called ``unsuccessful'', the iterate is rejected (\(\vek{x}_{k+1} = \vek{x}_k\)) and \(\sigma_k\) is increased.
This \textbf{a}daptively chosen \textbf{r}egularization parameter gives the AR\(p\) method its name \cite{birgin_worst-case_2017,cartis_evaluation_2022}.

The advantage of having access to higher derivatives can be quantified in two different ways on a theoretical level:
in terms of \emph{global} complexity guarantees, giving an upper bound on the number of iterations until an approximate stationary point is found, and in terms of \emph{local} convergence rates when the algorithm is already close to the solution.
For both cases it is standard to consider the class of \(p\)-times differentiable functions that have a Lipschitz continuous \(p\)th derivative and are lower bounded.
On these functions the \emph{global} complexity literature established that the AR\(p\) method has optimal worst-case global iteration complexity among all possible \(p\)th-order methods as it needs at most \(O(\e^{-(p+1)/p})\) iterations to find an approximate stationary point with a gradient norm bounded by \(\e > 0\)~\cite{birgin_worst-case_2017,cartis_sharp_2020,cartis_evaluation_2022}.
Clearly, this bound improves as \(p\) increases.

The \emph{local} convergence for these tensor methods has been analysed for the subclass of uniformly convex functions by Doikov and Nesterov~\cite{doikov_local_2022}.
A differentiable function is uniformly convex of degree \(q \geq 2\) if there exists a \(\mu > 0\) such that \(f(\vek{y}) \geq f(\vek{x}) + \nabla f(\vek{x})^T (\vek{y} - \vek{x}) + \mu \norm{\vek{y} - \vek{x}}^q\) for all \(\vek{x}, \, \vek{y} \in \R^n\).
This means that \(f\) is strictly convex and grows at least like \(O(\norm{\vek{x} - \vek{x}_*}^q)\) around any stationary point \(\vek{x}_*\).
Doikov and Nesterov~\cite{doikov_local_2022} additionally assume that the function's \(p\)th derivative is Lipschitz continuous with constant \(L_p\) and that the regularization parameter $\sigma_k$ is constant and sufficiently large (\(\sigma_k \geq \frac{p L_p}{(p+1)!}\)).
The latter condition allows them to analyse a non-adaptive version of AR\(p\) where \(\vek{y}_k\) is the global minimizer of the local model at each step, and it is always accepted.
For this method the authors obtain
\begin{subequations}
    \label{eqn:doikov-local-convergence}
    \begin{align}
    f(\vek{x}_{k+1}) - f^* &= O \Paren{(f(\vek{x}_k) - f^*)^{\frac{p}{q-1}}}, \label{eqn:doikov-local-convergence-function-value} \\
     \norm{\nabla f(\vek{x}_{k+1})} &= O \Paren{\norm{\nabla f(\vek{x}_k)}^{\frac{p}{q-1}}}. \label{eqn:doikov-local-convergence-gradient-norm}
    \end{align}
\end{subequations}
where \(f^*\) is the global minimum of \(f\).
\Cref{eqn:doikov-local-convergence-function-value,eqn:doikov-local-convergence-gradient-norm} imply that, if \(\vek{x}_k\) is close enough to the unique global minimizer \(\vek{x}_*\) and \(p > q-1\), then both the function value and the gradient norm converge superlinearly with order \(\frac{p}{q-1}\).

In this paper, we generalize the results of Doikov and Nesterov~\cite{doikov_local_2022} by
\begin{enumerate}
    \item considering the AR\(p\) algorithm, with its \emph{adaptively chosen regularization},
    \item requiring only \emph{approximate local model minimizers} at each iteration, and
    \item covering the case of \emph{locally uniformly convex functions}.
\end{enumerate}
In the local convergence setting where \(\vek{x}_0\) is sufficiently close to a local minimizer \(\vek{x}_*\), we obtain that the gradient norms converge with order~\(\frac{p}{q-1}\) when $p > q-1$ during successful iterations.
By additionally assuming which model minimizer \(\vek{y}_k\) is chosen when faced with ambiguity one can ensure that all iterations are successful, and we achieve the same local \(\frac{p}{q-1}\)th-order convergence of iterates, function values and gradient norms.

To put this into perspective, let us compare this convergence rate with the one of Newton's method.
It is well-known that Newton's method converges quadratically when the Hessian at the limit point is nonsingular or, equivalently, when \(f\) is locally uniformly convex of order \(q=2\).
In this case, the access to additional derivatives gives a speed-up from quadratic (Newton's method) to \(p\)th-order convergence (AR\(p\)).
Therefore, to achieve a tolerance of \(\e > 0\) on the norm of the gradient, both methods require \(O(\log(\log(\e^{-1})))\) iterations, but $p$th-order methods improve the coefficient of the  leading-order term from \(1/\log(2)\) to \(1/\log(p)\).
The advantage is even greater when the Hessian at \(\vek{x}_*\) is singular, as Newton's method only converges linearly in this case whereas the access to enough derivatives (\(p > q-1\)) allows superlinear convergence of order \(\frac{p}{q-1}\).
The maximum number of iterations are then bounded by \(O(\log(\e^{-1}))\) for Newton and \(O(\log(\log(\e^{-1})))\) for AR\(p\).

There are other works that cover the local convergence of regularization methods:
cubic regularization methods were introduced by Griewank \cite{griewank_modification_1981} and later independently revisited by Nesterov and Polyak~\cite{nesterov_cubic_2006}.
Cartis, Gould and Toint~\cite{cartis_adaptiveI_2011,cartis_adaptiveII_2011} added adaptivity of the regularization parameter and arrived at the ARC method (AR\(2\) in our notation).
Nesterov and Polyak~\cite[Theorem 3]{nesterov_cubic_2006} and Cartis, Gould and Toint~\cite[Corollary 4.10]{cartis_adaptiveI_2011} show quadratic local convergence of these second-order methods when the Hessian is positive definite at the minimizer, both when \(\sigma_k\) converges to zero and when it stays bounded away from zero.

Yue et al.~\cite{yue_quadratic_2019} go even further and are able to show that positive definiteness of the Hessian at the minimizer is not necessary for quadratic convergence of AR2, but rather that a local error bound condition (EB) suffices.
This agrees with earlier work by Fan and Yuan~\cite{fan_regularized_2014} on quadratic convergence of regularized second-order methods and work by Bellavia and Morini~\cite{bellavia_strong_2015} on quadratic convergence of adaptively regularized methods for nonlinear least squares.
Rebjock and Boumal~\cite{rebjock_fast_2024} prove the equivalence of different local properties including the error bound assumption.
With it, they conclude that in the case where the set of local minimizers with a certain function value is a submanifold \(\mathcal{S}\) and only the eigenvalues corresponding to normal directions are bounded away from zero, quadratic convergence of the AR2 method still holds.
In fact, for any neighbourhood \(\mathcal{U}\) of a local minimizer \(\bar{\vek{x}} \in \mathcal{S}\) there is a neighbourhood \(\mathcal{V}\) of \(\bar{\vek{x}}\) such that if an iterate enters \(\mathcal{V}\) then the sequence of AR2 iterates converges to some \(\vek{x}_* \in \mathcal{S} \cap \mathcal{U}\) Q-quadratically\footnote{Throughout the paper we use the terms Q-convergence and R-convergence as in \cite[section A.2]{nocedal_numerical_2006}. In particular, we say that a sequence converges linearly if it converges \emph{at least} linearly, without claiming that it does not converge at a higher-order rate. In this parlance, even a quadratically converging sequence also converges linearly.}.
To achieve this, the only assumption on how the regularization parameter \(\sigma_k\) is updated is that it does not increase during successful iterations and that it is always lower bounded by some \(\sigma_{\min} > 0\).

Cartis, Gould and Toint~\cite{cartis_evaluation_2022} prove complexity bounds for the AR\(p\) method on nonconvex objective functions under various assumptions.
In~\cite[Section 5.3]{cartis_evaluation_2022} they discuss measure-dominated problems.
As we will show in \cref{thm:uniform-convexity-properties}, functions that are uniformly convex of order \(q \geq 2\) are also gradient-dominated of level \(\frac{q}{q-1}\).
Using the fact that the function decrease of AR\(p\) methods can be lower bounded in a certain way, the results of \cite{cartis_evaluation_2022} imply that AR\(p\) methods are guaranteed to converge at least superlinearly if \(p > q-1\), at least linearly if \(p = q-1\) and at least sublinearly if \(p < q-1\).
More precisely, the corresponding upper bounds are \(O(\log(\log(\e^{-1})))\), \(O(\log(\e^{-1}))\) and \(O(\e^{-\frac{p+1}{p} + \frac{q}{q-1}})\) iterations, respectively, to achieve an error tolerance \(\norm{\nabla f(\vek{x})} < \e\) in the worst-case.
The \(p > q-1\) case matches the discussion above.
Note though, that by assuming different properties of $f$ (local uniform convexity instead of global gradient domination) and focusing on the gradient norm convergence, we can derive sharp bounds for the order of convergence missing from the results in \cite{cartis_evaluation_2022}.

The remainder of the paper is organized as follows: in \cref{sec:arp}, we describe the AR\(p\) algorithm and state its properties under basic assumptions.
Afterwards, an example, which motivates our paper, shows that even asymptotically, not all iterations need to be successful for higher-order methods unlike in the second-order case.
\Cref{sec:local-convergence} contains local convergence results for AR\(p\) under local uniform convexity assumptions on \(f\), both when not all iterations are successful and when this property is restored by additional assumptions on the choice of the model minimizer.
The sharpness of the derived convergence rates is examined in \cref{sec:sharpness} and the rates are illustrated with experiments in \cref{sec:numerical-illustration}.
We finally discuss the differences between our own analysis and other works in the literature and make some concluding remarks in \cref{sec:discussion}.

\section{\texorpdfstring{The AR\(p\) algorithm}{The ARp algorithm}}\label{sec:arp}
\begin{algorithm2e}
    \KwParameters{\(0 < \eta_1 \leq \eta_2 < 1\), \(0 < \gamma_1 \leq 1 < \gamma_2\), \(\theta > 0\)}
    \KwIn{\(\vek{x}_0 \in \R^n\), \(\sigma_0 > 0\)}
    \For{\(k = 0, 1, \dots\)}{
        \If{\(\nabla f(\vek{x}_k) = \vek{0}\)}{
            \Return{}
        }
        Find a local minimizer \(\vek{y}_k\) of \(m_k\) defined according to \cref{eqn:arp-model} that satisfies
            \begin{equation}\label{eqn:approximate-minimizer}
                m_k(\vek{y}_k) < m_k(\vek{x}_k) \quad \text{and} \quad \norm{\nabla m_k (\vek{x}_k)} \leq \theta \norm{\vek{x}_k - \vek{y}_k}^p
            \end{equation} \;
        \vspace{-\baselineskip}
        Set \(\rho_k = \dfrac{f(\vek{x}_k) - f(\vek{y}_k)}{t_k(\vek{x}_k) - t_k(\vek{y}_k)}\) \;
        \uIf{\(\rho_k \geq \eta_2\)}{
            Set \(\vek{x}_{k+1} = \vek{y}_k\) and \(\sigma_{k+1} = \gamma_1 \sigma_k\)  \tcp*[r]{very successful iteration}
        }
        \uElseIf{\(\rho_k \geq \eta_1\)}{
            Set \(\vek{x}_{k+1} = \vek{y}_k\) and \(\sigma_{k+1} = \sigma_k\)  \tcp*[r]{successful iteration}
        }
        \Else{
            Set \(\vek{x}_{k+1} = \vek{x}_k\) and \(\sigma_{k+1} = \gamma_2 \sigma_k\) \tcp*[r]{unsuccessful iteration}
        }
    }
	\caption{AR\(p\) algorithm}\label{alg:arp}
\end{algorithm2e}

In this paper, we aim to understand the asymptotic behaviour of the AR\(p\) algorithm, which is given in \cref{alg:arp}.
To this end, we assume that \cref{alg:arp} generates an infinite sequence \(\vek{x}_k\), i.e., that no iterate ever satisfies the termination condition \(\nabla f(\vek{x}_k) = \vek{0}\) exactly.

At iteration \(k\) the local Taylor expansion \(t_k\) and the local regularized model \(m_k\) at \(\vek{x}_k\) are defined as
\begin{subequations}
    \label{eqn:models}
    \begin{align}
        t_k(\vek{y}) &= t_{\vek{x}_k}^p(\vek{y}) = f(\vek{x}_k) + \sum_{j=1}^p \frac{1}{j!} \nabla^j f(\vek{x}_k) [\vek{y} - \vek{x}_k]^j \label{eqn:taylor-model} \\
        m_k(\vek{y}) &= m_{\vek{x}_k, \sigma_k}^p(\vek{y}) = t_{\vek{x}_k}^p(\vek{y}) + \sigma_k \norm{\vek{y} - \vek{x}_k}^{p+1} \label{eqn:arp-model}
    \end{align}
\end{subequations}
for \(\vek{y} \in \R^n\).
The tentative new iterate \(\vek{y}_k\) is computed such that it satisfies \cref{eqn:approximate-minimizer}.
We call any such point satisfying \cref{eqn:approximate-minimizer} an \emph{approximate minimizer} of \(m_k\).
After finding \(\vek{y}_k\), the algorithm checks whether \(\rho_k \geq \eta_1\) for this point,
i.e., whether the decrease in the objective function is at least a specified fraction \(\eta_1\) of the decrease predicted by the Taylor expansion \(t_k\).
The predicted decrease is always positive since \(t_k(\vek{y}_k) \leq m_k(\vek{y}_k) < m_k(\vek{x}_k) = t_k(\vek{x}_k)\).
If the sufficient decrease condition \(\rho_k \geq \eta\) is satisfied, the iteration is called ``successful'', and the iterate is accepted (\(\vek{x}_{k+1} = \vek{y}_k\)).
If moreover \(\rho_k \geq \eta_2\), then regularization parameter is decreased to \(\sigma_{k+1} = \gamma_1 \sigma_k\).
Otherwise, the iteration is called ``unsuccessful'', the iterate is rejected (\(\vek{x}_{k+1} = \vek{x}_k\)), and the regularization is increased to \(\sigma_{k+1} = \gamma_2 \sigma_k\).

In the remainder of this paper, we will always assume that \(\eta_1 = \eta_2 = \eta\) for convenience.
The presented convergence results also apply for the case when \(\eta_1 < \eta_2\), but the assumption simplifies the statements and proofs.
Compared to the AR\(p\) algorithms in \cite{birgin_worst-case_2017,cartis_sharp_2020,cartis_concise_2020} we have not included a \(\sigma_{\min}\) parameter that provides a lower bound for the regularization parameter.
Since we are analysing local convergence, there is no need to keep \(\sigma_k\) artificially away from zero.
As we will discuss in \cref{sec:right-minimizer}, under the right circumstances convergence at the optimal rate is guaranteed even when setting \(\sigma_k = 0\).

The regularization in \cref{eqn:arp-model} is chosen such that it compensates for the error in the Taylor approximation for the following function class:

\begin{definition}\label{def:pth-order-lipschitz}
    A function is \(p\)th-order Lipschitz continuous for some \(p \geq 2\) if it is \(p\) times continuous differentiable and there exists a constant \(L_p > 0\) such that
    \[
        \norm{\nabla^p f(\vek{x}) - \nabla^p f(\vek{y})} \leq L_p \norm{\vek{x} - \vek{y}}
    \]
    holds for any \(\vek{x}, \vek{y} \in \R^n\).
    For notational convenience, we define \(C^p_{\Lip}\) as the set of \(p\)th-order Lipschitz continuous functions.
\end{definition}

Under this assumption on \(f\), the error in the \(p\)th-order Taylor expansion and its derivatives can be bounded in the following way by Lemma 2.1 in \cite{cartis_sharp_2020}:
\begin{subequations}\label{eqn:taylor-error-bounds}
    \begin{align}
        \abs{f(\vek{y}) - t_{\vek{x}}^p(\vek{y})} &\leq \frac{L_p}{(p+1)!} \norm{\vek{x} - \vek{y}}^{p+1} \label{eqn:taylor-error-bound} \\
        \norm{\nabla f(\vek{y}) - \nabla t_{\vek{x}}^p(\vek{y})} &\leq \frac{L_p}{p!} \norm{\vek{x} - \vek{y}}^p \label{eqn:taylor-gradient-error-bound} \\
        \norm{\nabla^2 f(\vek{y}) - \nabla^2 t_{\vek{x}}^p(\vek{y})} &\leq \frac{L_p}{(p-1)!} \norm{\vek{x} - \vek{y}}^{p-1} \label{eqn:taylor-hessian-error-bound}
    \end{align}
\end{subequations}
In particular, when \(\sigma_k \geq L_p / (p+1)!\), then \(m_k(\vek{y}) \geq f(\vek{y})\) for all \(\vek{y} \in \R^n\).\footnote{Note that the denominators in \cref{eqn:taylor-error-bounds} differ from those found in the corresponding bounds in \cite{birgin_worst-case_2017,cartis_second-order_2018,cartis_concise_2020}, by a factor of \(p+1\), \(p\) and \(p-1\) respectively.
This is because the bounds derived in these papers are not sharp in terms of the constants.
Therefore, some fundamental lemmas in our paper differ in the constants involved compared to their counterparts in \cite{birgin_worst-case_2017,cartis_second-order_2018,cartis_concise_2020}.}

The following fundamental result shows that \(\sigma_k\) cannot become arbitrarily large throughout the course of the algorithm.

\begin{lemma}[Lemma 2.2 in \cite{birgin_worst-case_2017}, Lemma 3.2 in \cite{cartis_concise_2020}]\label{thm:sigma-max}
    For \(f \in C^p_{\Lip}\) the $k$th iteration is guaranteed to be successful (\(\rho_k \geq \eta\)) whenever \(\sigma_k \geq \frac{1}{1-\eta} \frac{L_p}{(p+1)!}\).
    As a consequence, the regularization parameter stays upper bounded:
    \[
        \sigma_k \leq \sigma_{\max} \deq \max \left\{\sigma_0, \ \frac{\gamma_2}{1-\eta} \frac{L_p}{(p+1)!} \right\}
    \]
    for all \(k \in \N\).
\end{lemma}

\begin{proof}
    It suffices to show that \(\abs{\rho_k - 1} \leq 1 - \eta\) whenever \(\sigma_k \geq \frac{1}{1-\eta} \frac{L_p}{(p+1)!}\).
    If \(\sigma_k\) is this large, we require
    \begin{align*}
        \abs{\rho_k - 1} &= \Abs{\frac{f(\vek{x}_k) - f(\vek{y}_k)}{t_k(\vek{x}_k) - t_k(\vek{y}_k)} - 1} = \frac{\abs{t_k(\vek{y}_k) - f(\vek{y}_k)}}{\abs{t_k(\vek{x}_k) - t_k(\vek{y}_k)}} \\
        &= \frac{\abs{t_k(\vek{y}_k) - f(\vek{y}_k)}}{\abs{m_k(\vek{x}_k) - m_k(\vek{y}_k) + \sigma_k \norm{\vek{y}_k - \vek{x}_k}^{p+1}}} \\
        &\leq \frac{\Paren{L_p / (p+1)!} \norm{\vek{y}_k - \vek{x}_k}^{p+1}}{\sigma_k \norm{\vek{y}_k - \vek{x}_k}^{p+1}} = \frac{L_p}{(p+1)!} \frac{1}{\sigma_k} \leq 1 - \eta
    \end{align*}
    using \(t_k(\vek{x}_k) = m_k(\vek{x}_k) = f(\vek{x}_k)\), \cref{eqn:arp-model,eqn:approximate-minimizer,eqn:taylor-error-bound}.
    The update mechanism of \cref{alg:arp} then ensures that the regularization parameter is never increased beyond \(\sigma_{\max}\).
\end{proof}

Note that if \(\gamma_1 < 1\), then starting at \(\sigma_0 \geq \frac{1}{1-\eta} \frac{L_p}{(p+1)!}\) the regularization parameter is decreased until it becomes smaller than this value.
Afterwards, it never increases beyond \(\frac{\gamma_2}{1-\eta} \frac{L_p}{(p+1)!}\) as shown above.
Therefore, asymptotically the regularization parameter can be bounded independently of \(\sigma_0\):
\[
  \limsup_{k \to \infty} \sigma_k \leq \frac{\gamma_2}{1-\eta} \frac{L_p}{(p+1)!}
\]

\subsection{Regularization parameter oscillations}\label{sec:example-sigma-oscillation}

In order to establish the difficulties inherent in proving local convergence of \(p\)th-order tensor methods for \(p \geq 3\) it is helpful to consider a simple one-dimensional example.
It shows that when \(\vek{x}_k\) is close to a non-degenerate global minimizer of the objective function, the \(p\)th-order Taylor expansion is locally strongly convex but does not have to be globally convex or even lower bounded, unlike for the second-order case.
This impacts the local convergence when \(\vek{y}_k\) is chosen as the global minimizer of \(m_k\).

\begin{figure}
    \centering
    \subimport{./plots}{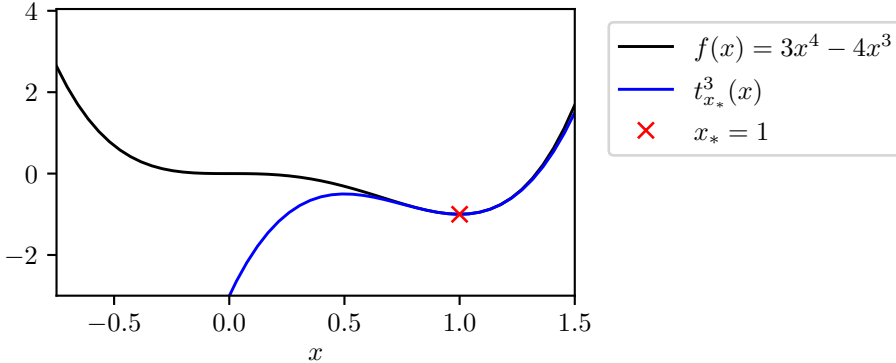}
    \vspace{-1em}
    \caption{Example function and its third-order Taylor expansion at the global minimizer \(x_* = 1\).}
    \label{fig:example1}
\end{figure}

\begin{example}\label{exm:sigma-oscillations}
    Consider the objective function \(f(x) = 3x^4 - 4x^3\) shown in \cref{fig:example1}.
    Its only stationary points are \(x = 0\) and \(x = 1\), of which the first is a saddle point and the second is the global minimizer \(x_*\).
    The objective function and its derivatives at \(x_*\) are \(f(x_*) = -1\), \(f'(x_*) = 0\), \(f''(x_*) = 12\), and \(f'''(x_*) = 48\).
    Therefore, \(x_*\) is a nondegenerate local minimizer in the sense that the Hessian is positive definite at \(x_*\) and the function is locally strongly convex around that point.
    Nonetheless, because \(f'''(x_*) \neq 0\) the third-order Taylor expansion at the minimizer is neither convex nor lower bounded.
    We will show that when \(y_k\) is chosen as the global minimizer of \(m_k\), then there exist parameters \(\sigma_0\), \(\eta\), \(\gamma_1\), and \(\gamma_2\) such that only every second iteration is successful for \emph{any} \(x_0\) close enough to \(x_*\).
    This is because on every other iteration, the regularization parameter is so small that the global minimizer of the model lies outside the convex region around the expansion point.

    Note that when \(\sigma = 3\) the local model is exact: \(m_{x, 3}^3(y) = t_{x}^3(y) + 3 (y - x)^4 = f(y).\)
    The last equality is easy to see by considering that both sides are quartic polynomials and comparing the first four derivatives at \(y = x\).
    It implies
    \[
        m_k(y) = f(y) + (\sigma_k - 3) (y - x)^4.
    \]
    For the given objective function we have \(L_3 = 4! \cdot 3\).
    By virtue of \cref{thm:sigma-max} we know that for \(\sigma_k \geq \frac{L_3}{4! (1 - \eta)}\) every iteration will be successful (\(\rho \geq \eta\)) regardless of which minimizer of \(m_k\) is chosen.
    Let \(\eta = 1/2\), then we obtain that the iteration is successful whenever \(\sigma_k \geq 6\).
    Next, we show that the global minimizer \(y_k\) of the model \(m_k\) leads to an unsuccessful iteration whenever \(\sigma_k \leq 2\) and \(x_k \in [\frac{4}{5}, \frac{6}{5}]\).
    In this case, we focus on the interval \([0, \frac{4}{3}]\) where \(f(x_k) < 0\). Notice that
    \[
        \min_{y \in [0, \frac{4}{3}]} m_k(y) = \min_{y \in [0, \frac{4}{3}]} f(y) + (\sigma_k - 3) (y - x_k)^4 \geq -1 + (\sigma_k - 3)\Paren{\frac{6}{5}}^4
    \]
    and
    \[
        m_k(-1) = f(-1) + (\sigma_k - 3)(x_k + 1)^4 \leq 7 + (\sigma_k - 3)\Paren{\frac{9}{5}}^4.
    \]
    The difference between these bounds is
    \[
        \Paren{\min_{y \in [0, \frac{4}{3}]} m_k(y)} - m_k(-1) \geq -8 + (\sigma_k - 3)\Paren{\Paren{\frac{6}{5}}^4 - \Paren{\frac{9}{5}}^4} > 0.
    \]
    This implies that any global minimizer \(y_k\) of \(m_k\) is outside the interval \([0, \frac{4}{3}]\), which means \(f(y_k) > 0 > f(x_k)\).
    Therefore, \(\rho_k < 0\) and the iteration is unsuccessful.

    Now, take \(\sigma_0 = 6\), \(\eta = \frac{1}{2}\), \(\gamma_1 = \frac{1}{3}\), \(\gamma_2 = 3\), and \(x_0\) with \(f(x_0) \leq -0.9\).
    Because \(\{x \in \R \mid f(x) \leq -0.9\} \subset [\frac{4}{5}, \frac{6}{5}]\) and iterates of \cref{alg:arp} have monotonically decreasing function values, all iterates lie in this interval.
    The first iteration with \(\sigma_0 = 6\) is successful, so the regularization parameter is decreased to \(\sigma_1 = 2\).
    The parameter is now too small (\(\sigma_1 \leq 2\)), the iteration is unsuccessful, and the parameter is increased to \(\sigma_2 = 6\).
    At this point the same cycle repeats ad infinitum.
    The points $x_k$ converge to $x_*$, but only every second iteration is successful and makes progress.
\end{example}

This example gives an important insight into the AR\(p\) method:
when choosing a global minimizer of \(m_k\) at each iteration, the regularization parameter \(\sigma_k\) will not in general go to zero asymptotically even when the objective function is locally strongly convex around the minimizer.
Since the global minimizer satisfies \cref{eqn:approximate-minimizer}, we conclude that it is impossible to show for the fully adaptive\footnote{By fully adaptive AR\(p\) we mean the case \(\gamma_1 < 1\), i.e., the regularization parameter is not just increased during unsuccessful iterations but also strictly decreased during successful iterations. In this terminology, the \(\gamma_1 = 1\) case might be called semi-adaptive.} AR\(p\) method in \cref{alg:arp} that asymptotically all iterations will be successful without specifying which minimizer of the local model needs to be chosen.

\section{Local convergence for locally uniformly convex functions}\label{sec:local-convergence}

We follow Doikov and Nesterov~\cite{doikov_local_2022} by considering the property of uniform convexity of order \(q\).
For our purposes though, we only require the function to be \emph{locally} uniformly convex in a closed convex neighbourhood around a local minimizer \(\vek{x}_*\):

\begin{definition}\label{def:uniform-convexity}
    A function is locally uniformly convex of order \(q \geq 2\) around a local minimizer \(\vek{x}_*\) if there exist parameters \(\mu_q, r_q > 0\) such that
    \begin{equation}\label{eqn:uniform-convexity-gradient}
        (\nabla f(\vek{x}) - \nabla f(\vek{y}))^T (\vek{x} - \vek{y}) \geq \mu_q \norm{\vek{x} - \vek{y}}^q \quad \forall \vek{x}, \vek{y} \in \ball(\vek{x}_*, r_q),
    \end{equation}
    where \(\ball(\vek{x}_*, r_q)\) is a closed ball of radius \(r_q\) around \(\vek{x}_*\).
    For notational convenience, we write \(M^q(f)\) for the set of local minimizers \(\vek{x}_*\) such that \cref{eqn:uniform-convexity-gradient} holds.
\end{definition}

Note that by virtue of \cite[Lemma 4.2.1]{nesterov_lectures_2018}, \cref{eqn:uniform-convexity-gradient} implies
\begin{equation}\label{eqn:uniform-convexity-function}
    f(\vek{y}) \geq f(\vek{x}) + \nabla f(\vek{x})^T(\vek{y} - \vek{x}) + \frac{\mu_q}{q} \norm{\vek{y} - \vek{x}}^q  \quad \forall \vek{x}, \vek{y} \in \ball(\vek{x}_*, r_q),
\end{equation}
and simultaneously \cref{eqn:uniform-convexity-function} implies \cref{eqn:uniform-convexity-gradient} (albeit with a different constant \(\mu_q\)). Therefore, the two inequalities define the same notion of uniform convexity.

Local uniform convexity of order \(q\) implies local strict convexity and therefore that \(\vek{x}_*\) is an isolated minimizer, but for \(q > 2\) the assumption is weaker than local strong convexity.
Indeed, local strong convexity, namely \(\nabla^2 f(\vek{x}) \succeq \mu_2 \mat{I}\) for all \(\vek{x} \in \ball(\vek{x}_*, r_q)\), is equivalent to \cref{eqn:uniform-convexity-gradient} with \(q = 2\).
By analysing uniformly convex functions instead of strongly convex functions we are able to show superlinear convergence even for degenerate minimizers in certain cases.

\begin{lemma}\label{thm:uniform-convexity-properties}
    Let \(\vek{x}_* \in M^q(f)\).
    Then \(f\) satisfies
    \begin{align}
        f(\vek{x}) - f(\vek{x}_*) &\geq \frac{\mu_q}{q} \norm{\vek{x} - \vek{x}_*}^q, \label{eqn:uniform-convexity-step-bounded-by-function} \\
        \norm{\nabla f(\vek{x})} &\geq \mu_q \norm{\vek{x} - \vek{x}_*}^{q-1}, \text{and} \label{eqn:uniform-convexity-step-bounded-by-gradient} \\
        f(\vek{x}) - f(\vek{x}_*) &\leq \frac{q-1}{q} \Paren{\frac{1}{\mu_q}}^{\frac{1}{q-1}} \norm{\nabla f(\vek{x})}^{\frac{q}{q-1}} \label{eqn:uniform-convexity-function-bounded-by-gradient}
    \end{align}
    for all \(\vek{x} \in \ball(\vek{x}_*, r_q)\).
    Moreover, for \(f \in \mathcal{C}^2\) there exists some \(\nu > 0\) such that
    \begin{align}
        f(\vek{x}) - f(\vek{x}_*) &\leq \frac{\nu}{2} \norm{\vek{x} - \vek{x}_*}^2 \label{eqn:bounded-hessian-function-bound} \\
        \norm{\nabla f(\vek{x})} &\leq \nu \norm{\vek{x} - \vek{x}_*} \label{eqn:bounded-hessian-gradient-bound}
    \end{align}
    hold for all \(\vek{x} \in \ball(\vek{x}_*, r_q)\).
\end{lemma}

\begin{proof}
    The first inequality follows from \cref{eqn:uniform-convexity-function} by substituting \(\vek{x} = \vek{x}_*\) and \(\vek{y} = \vek{x}\).
    The second follows from \cref{eqn:uniform-convexity-gradient} by considering
    \[
        \norm{\nabla f(\vek{x})} \norm{\vek{x} - \vek{x}_*} \geq \Paren{\nabla f(\vek{x}) - \nabla f(\vek{x}_*)}^T (\vek{x} - \vek{x}_*) \geq \mu_q \norm{\vek{x} - \vek{x}_*}^q.
    \]
    By \cite[Lemma 4.2.2]{nesterov_lectures_2018}, we obtain \cref{eqn:uniform-convexity-function-bounded-by-gradient}.

    For the second part, the 2-norm of \(\nabla^2 f\) is bounded by some \(\nu \geq 0\) on the compact set \(\ball(\vek{x}_*, r_q)\).
    Equivalently, \(\nabla f\) is locally Lipschitz continuous with \(L_1 = \nu\) on that ball.
    Apply \cref{eqn:taylor-error-bound,eqn:taylor-gradient-error-bound} to the Taylor expansion at \(\vek{x}_*\) to obtain \cref{eqn:bounded-hessian-function-bound,eqn:bounded-hessian-gradient-bound}.
\end{proof}

\subsection{Convergence without assumptions on the model minimizer}

We first consider the best possible convergence results that can be derived for \cref{alg:arp} as stated.
In particular, there are no further restrictions on \(\vek{y}_k\) other than \cref{eqn:approximate-minimizer}.
Since the Taylor expansions are not necessarily convex, even when \(\vek{x}_k\) is close to \(\vek{x}_*\), there can be many local minimizers of \(m_k\) as long as \(\sigma_k\) is small enough.

The common assumptions for all the results in this section and the next are that \(f\) is \(p\)th-order Lipschitz continuous (abbreviated as \(f \in C^p_{\Lip}\), see \cref{def:pth-order-lipschitz}), that \(f\) is locally uniformly convex of order \(q\) around the local minimizer \(\vek{x}_*\) (abbreviated as \(\vek{x}_* \in M^q(f)\), see \cref{def:uniform-convexity}) and that \(p > q-1\).
These will be stated explicitly at the start of each statement.
Additionally, we silently assume that \(\vek{x}_k\) and \(\vek{y}_k\) are generated by \cref{alg:arp} using all \(p\) derivatives of the objective function and that the algorithm never terminates, i.e. that \(\nabla f(\vek{x}_k) \neq \vek{0}\) for all \(k\).

The following two lemmas form the basis of our convergence analysis.
The first result bounds the gradient norm at the tentative new iterate \(\vek{y}_k\) in terms of the gradient norm at the current iterate \(\vek{x}_k\) during successful iterations, while the second provides a corresponding bound for the function value gap \(f(\vek{x}_k) - f(\vek{x}_*)\).

\begin{lemma}\label{thm:gradient-convergence-pq-rate}
    Let \(f \in C^p_{\Lip}\) and \(\vek{x}_* \in M^q(f)\).
    Assume that \(\vek{x}_k, \vek{y}_k \in \ball(\vek{x}_*, r_q)\) and \(f(\vek{y}_k) \leq f(\vek{x}_k)\), then the gradients satisfy
    \begin{align}
        \norm{\nabla f(\vek{y}_k)} &\leq R_k \norm{\vek{y}_k - \vek{x}_k}^p \leq R_k \Paren{\frac{q}{\mu_q}}^{\frac{p}{q-1}} \norm{\nabla f(\vek{x}_k)}^{\frac{p}{q-1}} \label{eqn:gradient-yk-convergence-pq-rate}
    \end{align}
    where \(R_k \deq \paren[\big]{L_p / p! + \theta + (p+1) \sigma_k}\).
    If all iterates stay within \(\ball(\vek{x}_*, r_q)\), then
    \begin{equation}\label{eqn:gradient-convergence-pq-rate}
        \norm{\nabla f(\vek{x}_{k+1})} \leq R_{\max} \Paren{\frac{q}{\mu_q}}^{\frac{p}{q-1}} \norm{\nabla f(\vek{x}_k)}^{\frac{p}{q-1}}
    \end{equation}
    holds at every successful iteration for \(R_{\max} \deq \paren[\big]{L_p / p! + \theta + (p+1) \sigma_{\max}}\).
\end{lemma}

\begin{proof}
    Following the argument in \cite[Lemma 2.3]{birgin_worst-case_2017}, we have
    \begin{align*}
        \norm{\nabla f(\vek{y}_k)} &\leq \norm{\nabla f(\vek{y}_k) - \nabla t_k(\vek{y}_k)} + \norm{\nabla t_k(\vek{y}_k)} \\
        &\leq \norm{\nabla f(\vek{y}_k) - \nabla t_k(\vek{y}_k)} + \norm{\nabla m_k(\vek{y}_k)} + (p+1) \sigma_k \norm{\vek{y}_k - \vek{x}_k}^p \\
        &\leq \frac{L_p}{p!} \norm{\vek{y}_k - \vek{x}_k}^p + \theta \norm{\vek{y}_k - \vek{x}_k}^p + (p+1) \sigma_k \norm{\vek{y}_k - \vek{x}_k}^p \\
        &= \paren[\big]{L_p / p! + \theta + (p+1) \sigma_k} \norm{\vek{y}_k - \vek{x}_k}^p
    \end{align*}
    where we have used the triangle inequality, \cref{eqn:taylor-gradient-error-bound}, and \cref{eqn:approximate-minimizer}.

    Rearranging \cref{eqn:uniform-convexity-function} with \(\vek{x} = \vek{x}_k\) and \(\vek{y} = \vek{y}_k\) gives
    \[
        \nabla f(\vek{x}_k)^T (\vek{x}_k - \vek{y}_k) \geq f(\vek{x}_k) - f(\vek{y}_k) + \frac{\mu_q}{q} \norm{\vek{x}_k - \vek{y}_k}^q \geq \frac{\mu_q}{q} \norm{\vek{x}_k - \vek{y}_k}^q
    \]
    and therefore \(\norm{\nabla f(\vek{x}_k)} \geq \frac{\mu_q}{q} \norm{\vek{x}_k - \vek{y}_k}^{q-1}\) by \(f(\vek{y}_k) \leq f(\vek{x}_k)\).
    Combine this result with the bound on \(\norm{\nabla f(\vek{y}_k)}\) to obtain \cref{eqn:gradient-yk-convergence-pq-rate}.

    By \cref{thm:sigma-max} we know \(\sigma_k \leq \sigma_{\max}\), and during successful iterations \(\vek{x}_{k+1} = \vek{y}_k\) and \(f(\vek{x}_{k}) - f(\vek{y}_k) \geq \eta \Paren{t_k(\vek{x}_k) - t_k(\vek{y}_k)} \geq \eta \Paren{m_k(\vek{x}_k) - m_k(\vek{y}_k)} > 0 \), so \cref{eqn:gradient-convergence-pq-rate} follows immediately.
\end{proof}

\begin{lemma}\label{thm:function-convergence-pq-rate}
    Let \(f \in C^p_{\Lip}\), \(\vek{x}_* \in M^q(f)\) and \(p > q-1\).
    Assume that \(\vek{x}_k, \vek{y}_k \in \ball(\vek{x}_*, r_q)\), \(f(\vek{y}_k) \leq f(\vek{x}_k)\) and that the norm of the gradient at \(\vek{x}_k\) is small enough.
    Then the function values satisfy
    \begin{equation}\label{eqn:function-yk-convergence-pq-rate}
        f(\vek{y}_k) - f(\vek{x}_*) \leq \frac{q-1}{q} \Paren{\frac{1}{\mu_q}}^{\frac{1}{q-1}} \Paren{\frac{2q}{\mu_q}}^{\frac{p}{q-1}} R_k^{\frac{q}{q-1}} \Paren{f(\vek{x}_k) - f(\vek{x}_*)}^{\frac{p}{q-1}}.
    \end{equation}

    In particular, if all iterates stay within \(\ball(\vek{x}_*, r_q)\) and \(\norm{\nabla f(\vek{x}_k)} \to 0\), then
    \begin{equation}\label{eqn:function-convergence-pq-rate}
        f(\vek{x}_{k+1}) - f(\vek{x}_*) \leq \frac{q-1}{q} \Paren{\frac{1}{\mu_q}}^{\frac{1}{q-1}} \Paren{\frac{2q}{\mu_q}}^{\frac{p}{q-1}} R_{\max}^{\frac{q}{q-1}} \Paren{f(\vek{x}_k) - f(\vek{x}_*)}^{\frac{p}{q-1}}
    \end{equation}
    holds at every successful iteration \(k\) for \(k\) large enough.
\end{lemma}
\begin{proof}
    Whenever the norm of the gradient at \(\vek{y}_k\) is small enough such that
    \begin{equation}\label{eqn:gradient-norm-condition-for-function-pq-convergence}
        \norm{\nabla f(\vek{y}_k)}^{\frac{p - q + 1}{p}} \leq \frac{\mu_q}{2q} \Paren{\frac{1}{R_{\max}}}^{\frac{q-1}{p}},
    \end{equation}
    then we can use \cref{eqn:uniform-convexity-function} and the first inequality in \cref{eqn:gradient-yk-convergence-pq-rate} to obtain
    \begin{align*}
        f(\vek{x}_k) - f(\vek{x}_*) &\geq f(\vek{x}_k) - f(\vek{y}_k) \\
        &\geq -\norm{\nabla f(\vek{y}_k)}\norm{\vek{x}_k - \vek{y}_k} + \frac{\mu_q}{q} \norm{\vek{x}_k - \vek{y}_k}^q \\
        &\geq \Paren{\frac{\mu_q}{q} \Paren{\frac{\norm{\nabla f(\vek{y}_k)}}{R_k}}^{\frac{q-1}{p}} - \norm{\nabla f(\vek{y}_k)}} \norm{\vek{x}_k - \vek{y}_k} \\
        &\geq \frac{\mu_q}{2q} \Paren{\frac{1}{R_k}}^{\frac{q}{p}} \norm{\nabla f(\vek{y}_k)}^{\frac{q}{p}}.
    \end{align*}
    \Cref{eqn:function-yk-convergence-pq-rate} follows from this and \cref{eqn:uniform-convexity-function-bounded-by-gradient} applied to \(\vek{y}_k\).
    Note that from \cref{eqn:gradient-convergence-pq-rate} we know that \cref{eqn:gradient-norm-condition-for-function-pq-convergence} holds as long as \(\norm{\nabla f(\vek{x}_k)}\) is small enough.

    As before \(\sigma_k \leq \sigma_{\max}\), the function value decreases in successful iterations and by assumption the gradient norms converge to zero, which means that the second claim follows from \cref{eqn:function-yk-convergence-pq-rate}.
\end{proof}

\Cref{thm:gradient-convergence-pq-rate,thm:function-convergence-pq-rate} show that asymptotically the norm of the gradients and the function value gaps converge at a Q-\(\frac{p}{q-1}\)th-order rate during successful iterations when \(p > q-1\).
Using local uniform convexity of \(f\), the convergence of the iterates is also of order \(\frac{p}{q-1}\) as shown by the following theorem.
The necessity of the assumption that all minimizers stay in the uniform convexity region around \(\vek{x}_*\) is clarified in \cref{thm:non-escape} and the preceding discussion.

\begin{theorem}\label{thm:convergence-pq-rate-successful}
    Let \(f \in C^p_{\Lip}\), \(\vek{x}_* \in M^q(f)\) and \(p > q-1\).
    Assume that all iterates stay in \(\ball(\vek{x}_*, r_q)\) and that \(\vek{x}_0\) is close enough to \(\vek{x}_*\). Then
    \begin{align*}
        \norm{\nabla f(\vek{x}_{k_i})} &\to 0 && \text{at a Q-} \tfrac{p}{q-1} \text{th-order rate} \\
        f(\vek{x}_{k_i}) &\to f(\vek{x}_*) && \text{at a Q-} \tfrac{p}{q-1} \text{th-order rate} \\
        \vek{x}_{k_i} &\to \vek{x}_* && \text{at an R-} \tfrac{p}{q-1} \text{th-order rate}
    \end{align*}
    as \(i \to \infty\), where \(\{k_1, k_2, \dots\}\) are the successful iterations.
\end{theorem}

\begin{proof}
    Note that the iterates stay constant in unsuccessful iterations, so it suffices to analyse how gradient norms, function values and iterates change from \(\vek{x}_k\) to \(\vek{x}_{k+1}\) in successful iterations.

    Since \(\nabla f\) is continuous, we know that for \(\vek{x}_0\) close enough to \(\vek{x}_*\), \(\norm{\nabla f(\vek{x}_0)}\) becomes small enough such that \cref{eqn:gradient-convergence-pq-rate} ensures at least linear decrease to zero at every successful iteration.
    This shows that \(\norm{\nabla f(\vek{x}_{k_i})}\) converges to zero, which in turn means that \cref{eqn:gradient-convergence-pq-rate} also implies the Q-superlinear convergence rate and that \cref{eqn:function-convergence-pq-rate} shows the same Q-superlinear rate for the function value gaps \(f(\vek{x}_{k_i}) - f(\vek{x}_*)\).

    The R-convergence of \(\vek{x}_{k_i}\) follows by \cref{thm:uniform-convexity-properties}, because
    \[
        \norm{\vek{x}_{k_i} - \vek{x}_*} \leq \mu_q^{-\frac{1}{q-1}} \norm{\nabla f(\vek{x}_{k_i})}^{\frac{1}{q-1}}.
    \]
    Even though the bound includes an exponent on \(\norm{\nabla f(\vek{x}_{k_i})}\) this does not change the order of convergence.
\end{proof}

Now that we know how fast gradient norms, function values and iterates converge during successful iterates, we can give a lower bound on the convergence rate of all iterates by exploiting that we can bound the number of unsuccessful iterations in terms of the number of successful ones.
The ratio of unsuccessful to successful iterations depends on the choice of \(\gamma_1\) and \(\gamma_2\).
In particular, when \(\gamma_1 = 1\) (\(\sigma_k\) is never decreased) then asymptotically all iterations are successful and we recover the full \(\frac{p}{q-1}\)th-order rates from \cref{thm:convergence-pq-rate-successful}.

\begin{theorem}\label{thm:convergence-pq-rate}
    Let \(f \in C^p_{\Lip}\), \(\vek{x}_* \in M^q(f)\) and \(p > q-1\).
    Assume that all iterates stay in \(\ball(\vek{x}_*, r_q)\) and that \(\vek{x}_0\) is close enough to \(\vek{x}_*\). Then
    \begin{align*}
        \norm{\nabla f(\vek{x}_k)} &\to 0 && \text{at an R-} \sqrt[\alpha+1]{\tfrac{p}{q-1}} \text{th-order rate} \\
        f(\vek{x}_k) &\to f(\vek{x}_*) && \text{at an R-} \sqrt[\alpha+1]{\tfrac{p}{q-1}} \text{th-order rate} \\
        \vek{x}_k &\to \vek{x}_* && \text{at an R-} \sqrt[\alpha+1]{\tfrac{p}{q-1}} \text{th-order rate}
    \end{align*}
    as \(k \to \infty\), where \(\gamma_1 = \gamma_2^{-\alpha}\).
\end{theorem}

\begin{proof}
    The assumptions are the same as those of \cref{thm:convergence-pq-rate-successful}, so we know the convergence rates when considering only the successful iterations.
    Let \(s_k\) and \(u_k\) be the number of successful and unsuccessful iterations from \(0\) to \(k-1\) respectively and let \(\{k_1, k_2, \dots\}\) be the indices of the successful iterations.
    Note that since the iterates do not change in unsuccessful iterations, we have that \(\vek{x}_k = \vek{x}_{k_i}\) for \(i = s_k\).

    Applying the same reasoning as \cite[Lemma 2.1]{cartis_adaptiveII_2011} and \cite[Lemma 2.4]{birgin_worst-case_2017}, we see that
    \[
        \sigma_0 \gamma_1^{s_k} \gamma_2^{u_k} = \sigma_k \leq \sigma_{\max}.
    \]
    Therefore, substituting \(\gamma_1 = \gamma_2^{-\alpha}\), taking logs on both sides and rearranging we get that \(u_k - \alpha s_k\) is upper bounded by the constant \(\log(\sigma_{\max} / \sigma_0) / \log(\gamma_2)\).
    Expressed in another way, using \(s_k + u_k = k\), we get
    \[
        s_k \geq \frac{k}{\alpha + 1} - \frac{\log(\sigma_{\max} / \sigma_0)}{(\alpha + 1) \log(\gamma_2)}.
    \]

    To establish \(r\)th-order rates of R-convergence for a sequence of values \((\e_k)_{k \in \N}\) that converge to zero, we need to show that asymptotically \(\log(\e_k^{-1})\) grows at least as fast as \(\Theta(r^k)\).
    By \cref{thm:convergence-pq-rate-successful} we already know that there exist \(c > 0\) such that
    \[
        \log(\norm{\nabla f(\vek{x}_{k_i})}^{-1}) \geq c \Paren{\tfrac{p}{q-1}}^i
    \]
    holds for all \(i\) large enough.
    This implies
    \begin{align*}
        \log(\norm{\nabla f(\vek{x}_k)}^{-1}) &= \log(\norm{\nabla f(\vek{x}_{k_{s_k}})}^{-1}) \geq c \Paren{\tfrac{p}{q-1}}^{s_k} \\
        &\geq c \Paren{\tfrac{p}{q-1}}^{\frac{k}{\alpha + 1} - \frac{\log(\sigma_{\max} / \sigma_0)}{(\alpha + 1) \log(\gamma_2)}} = c \Paren{\tfrac{p}{q-1}}^{- \frac{\log(\sigma_{\max} / \sigma_0)}{(\alpha + 1) \log(\gamma_2)}} \Paren{\sqrt[\alpha + 1]{\tfrac{p}{q-1}}}^k
    \end{align*}
    for all \(k\) large enough, which means \(\norm{\nabla f(\vek{x}_k)} \to 0\) at the claimed rate.
    Following the same arguments, we also get the correct convergence rates for \(f(\vek{x}_k) \to f(\vek{x}_*)\) and \(\vek{x}_k \to \vek{x}_*\) from \cref{thm:convergence-pq-rate-successful}.
\end{proof}

Note that even though the order of convergence decreased between \cref{thm:convergence-pq-rate-successful,thm:convergence-pq-rate}, it is still greater than one, since \(\alpha \geq 0\) and \(p > q-1\).
Therefore, the convergence stays superlinear.

A common assumption in the previous results is that all iterates stay in the neighbourhood of uniform convexity around \(\vek{x}_*\).
Clearly, this is unavoidable as we made very few assumptions about \(f\) outside this region.
In particular, as \(\vek{x}_*\) is only a local minimizer, a tentative iterate \(\vek{y}_k\) outside the neighbourhood might achieve an even lower function value than \(f(\vek{x}_*)\) and be accepted as \(\vek{x}_{k+1}\).
In such case, the algorithm does not converge to \(\vek{x}_*\).
Indeed, it will never return to the neighbourhood of \(\vek{x}_*\), since \(f(\vek{x}_k)\) is monotonically decreasing.
The next result shows that, if \(\vek{x}_0\) is close enough to \(\vek{x}_*\), then the case described above is the only case in which the algorithm leaves the neighbourhood of \(\vek{x}_*\).

\begin{theorem}\label{thm:non-escape}
    Let \(f \in C^p_{\Lip}\), \(\vek{x}_* \in M^q(f)\) and \(p > q-1\).
    If \(\vek{x}_0\) is close enough to \(\vek{x}_*\) (depending on \(\sigma_0\)), then either all iterates stay inside \(\ball(\vek{x}_*, r_q)\) or the first iterate outside this ball satisfies \(f(\vek{x}_k) < f(\vek{x}_*)\).
\end{theorem}

\begin{proof}
    The iterates only change when the iteration is successful, and it suffices to consider successful iterations. Let \(\{k_1, k_2, \dots\}\) be the set of iterations \(k\) where \(\rho_{k} \geq \eta\), then \(\vek{x}_{k_{i+1}} = \vek{x}_{k_i + 1}\).
    Moreover, let \(I \in \N_{> 0}\) be the first iteration such that \(f(\vek{x}_{k_I}) < f(\vek{x}_*)\) or infinity if no such iteration exists.
    We have that
    \begin{align*}
        f(\vek{x}_{k_i}) - f(\vek{x}_*) &\geq f(\vek{x}_{k_i}) - f(\vek{x}_{k_i + 1}) \geq \eta \Paren{t_{k_i}(\vek{x}_{k_i}) - t_{k_i}(\vek{x}_{k_i + 1})} \\
        &\geq \eta \sigma_{k_i} \norm{\vek{x}_{k_i} - \vek{x}_{k_i + 1}}^{p+1} \geq \eta \sigma_0 \gamma_1^i \norm{\vek{x}_{k_i} - \vek{x}_{k_i + 1}}^{p+1}
    \end{align*}
    for all \(i + 1 < I\) using \(\rho_{k_i} \geq \eta\), \cref{eqn:approximate-minimizer}, and the fact that \(\sigma_k\) is only decreased on successful iterations.
    Therefore,
    \begin{subequations}\label{eqn:distance-bound-successful-iteration}
        \begin{align}
            &\phantom{{}\leq{}} \norm{\vek{x}_{k_{i+1}} - \vek{x}_*} =  \norm{\vek{x}_{k_i + 1} - \vek{x}_*} \leq \norm{\vek{x}_{k_i + 1} - \vek{x}_{k_i}} + \norm{\vek{x}_{k_i} - \vek{x}_*} \\
            &\leq \Paren{\frac{f(\vek{x}_{k_i}) - f(\vek{x}_*)}{\eta \sigma_0 \gamma_1^i}}^{\frac{1}{p+1}} + \norm{\vek{x}_{k_i} - \vek{x}_*} \\
            &\leq \Paren{\frac{\frac{q-1}{q} \Paren{\frac{1}{\mu_q}}^{\frac{1}{q-1}}}{\eta \sigma_0} \frac{\norm{\nabla f(\vek{x}_{k_i})}^{\frac{q}{q-1}}}{\gamma_1^i}}^{\frac{1}{p+1}} + \Paren{\frac{1}{\mu_q}}^{\frac{1}{q-1}} \norm{\nabla f(\vek{x}_{k_i})}^{\frac{1}{q-1}},
        \end{align}
    \end{subequations}
    where the last inequality follows from \cref{thm:uniform-convexity-properties} assuming \(\vek{x}_{k_i}\) is inside \(\ball(\vek{x}_*, r_q)\).

    At this point, we cannot yet assert that either \cref{eqn:gradient-convergence-pq-rate} or \cref{eqn:distance-bound-successful-iteration} holds, since the former requires \(\vek{x}_{k_i + 1}, \vek{x}_{k_i} \in \ball(\vek{x}_*, r_q)\) and the latter \(\vek{x}_{k_i} \in \ball(\vek{x}_*, r_q)\).
    Still, choosing \(\vek{x}_0\) close enough to \(\vek{x}_*\) we have that \(\norm{\nabla f(\vek{x}_0)}\) gets arbitrarily small.
    We can select \(\vek{x}_0\) in such a way that it simultaneously satisfies
    \begin{enumerate}
        \item \(\vek{x}_0 \in \ball(\vek{x}_*, r_q)\),
        \item the right-hand side of \cref{eqn:distance-bound-successful-iteration} is smaller than \(r_q\) for \(i = 0\), and
        \item the right-hand side of \cref{eqn:gradient-convergence-pq-rate} is bounded by \(\gamma_1^{(q-1)/q} \norm{\nabla f(\vek{x}_0)}\) for \(k = 0\).
    \end{enumerate}
    By choosing \(\vek{x}_{k_0} = \vek{x}_0\) in such a way, it follows from \cref{eqn:distance-bound-successful-iteration} that \(\vek{x}_{k_1}\) lies inside \(\ball(\vek{x}_*, r_q)\) and therefore its gradient norm is bounded by \cref{eqn:gradient-convergence-pq-rate}.
    This implies that the gradient norm decreases at least by a factor of \(\gamma_1^{(q-1)/q}\).
    Then the right-hand side of \cref{eqn:distance-bound-successful-iteration} for \(i=1\) is smaller than the bound for \(i = 0\), which was already smaller than \(r_q\), and so \(\vek{x}_{k_2} \in \ball(\vek{x}_*, r_q)\).
    By induction, all successful iterates \(\vek{x}_{k_i}\) with \(i < I\) lie inside a circle of radius \(r_q\) around \(\vek{x}_*\).
    This proves the claim.
\end{proof}

\subsection{Choosing the right minimizer}\label{sec:right-minimizer}
The results of the previous subsection capture the worst-case convergence rate when all iterates stay in the region of uniform convexity.
They explicitly deal with the fact that even asymptotically not all iterations need to be successful.
However, this can be amended by specifying which local minimizer of \(m_k\) should be chosen, if there are multiple options.
In the following, we show that there always exists a local minimizer in the vicinity of \(\vek{x}_*\) that leads to a successful iteration (\cref{thm:close-minimizer,thm:connected-component-minimizer}).
Recall that in this paper an approximate model minimizer is simply one that satisfies \cref{eqn:approximate-minimizer}.
By selecting any such minimizer close to \(\vek{x}_*\), we can restore Q-\(\frac{p}{q-1}\)th-order convergence of the gradient norms regardless of the values of \(\gamma_1\) and \(\gamma_2\) (\cref{thm:convergence-pq-rate-good-min}).
As a stepping stone, we use that even though \(t_k\) is not necessarily convex around \(\vek{x}_*\), even when the expansion point is close to \(\vek{x}_*\),\footnote{As a counterexample consider $f(x) = x^4 + x^5$ (which is locally uniformly convex of order $q=4$ around $x_* = 0$) and its fourth-order Taylor expansion $t_k$ at a point $x_k < 0$. Analytically, we can deduce that $t_k''(0) = 20 x_k^3 < 0$ which means $t_k$ is not locally convex around $x_*$.} it does satisfy \cref{eqn:uniform-convexity-function} if one of \(\vek{x}\) and \(\vek{y}\) is the expansion point.

\begin{lemma} \label{thm:taylor-partial-uniform-convexity}
    Let \(f \in C^p_{\Lip}\), \(\vek{x}_* \in M^q(f)\) and \(p > q-1\).
    There exist \(\mu_q' < \mu_q\) and \(r_q' \leq r_q\) such that
    \begin{align}
        t_{\vek{x}}^p(\vek{x}) &\geq t_{\vek{x}}^p(\vek{y}) + \nabla t_{\vek{x}}^p(\vek{y})^T (\vek{x} - \vek{y}) + \frac{\mu_q'}{q} \norm{\vek{x} - \vek{y}}^q \label{eqn:taylor-uniform-convexity-x} \\
        t_{\vek{x}}^p(\vek{y}) &\geq t_{\vek{x}}^p(\vek{x}) + \nabla t_{\vek{x}}^p(\vek{x})^T (\vek{y} - \vek{x}) + \frac{\mu_q'}{q} \norm{\vek{y} - \vek{x}}^q \label{eqn:taylor-uniform-convexity-y}
    \end{align}
    hold for all \(\vek{x}, \vek{y} \in \ball(\vek{x}_*, r_q')\).
\end{lemma}

\begin{proof}
    Let \(\vek{x}, \vek{y} \in \ball(\vek{x}_*, r_q')\) for \(r_q' \leq r_q\) small enough such that
    \[
        \norm{\vek{x} - \vek{y}} \leq \norm{\vek{x} - \vek{x}_*} + \norm{\vek{y} - \vek{x}_*} \leq \sqrt[p-q+1]{\frac{\mu_q}{2q} \Big/ \Paren{\frac{L_p}{(p+1)!} + \frac{L_p}{p!}}}.
    \]
    In the first case we use \cref{eqn:taylor-error-bound}, \cref{eqn:taylor-gradient-error-bound}, and \cref{eqn:uniform-convexity-function} to obtain
    \[
        t_{\vek{x}}^p(\vek{x}) - t_{\vek{x}}^p(\vek{y}) - \nabla t_{\vek{x}}^p(\vek{y})^T (\vek{x} - \vek{y}) \geq \frac{\mu_q}{q} \norm{\vek{x} - \vek{y}}^q - \Paren{\frac{L_p}{(p+1)!} + \frac{L_p}{p!}}\norm{\vek{x} - \vek{y}}^{p+1}
    \]
    and similarly
    \[
        t_{\vek{x}}^p(\vek{y}) - t_{\vek{x}}^p(\vek{x}) - \nabla t_{\vek{x}}^p(\vek{x})^T (\vek{y} - \vek{x}) \geq \frac{\mu_q}{q} \norm{\vek{x} - \vek{y}}^q - \frac{L_p}{(p+1)!}\norm{\vek{x} - \vek{y}}^{p+1}.
    \]
    With above bound on \(\norm{\vek{x}- \vek{y}}\) we conclude that \cref{eqn:taylor-uniform-convexity-x,eqn:taylor-uniform-convexity-y} hold for \(\mu_q' = \mu_q / 2\).
\end{proof}

\begin{lemma}\label{thm:boundary-sublevel-set}
    Let \(f \colon \R^n \to \R\) be continuous and \(\mathcal{L}_f = \{\vek{y} \in \R^n \mid f(\vek{y}) \leq c\}\) a sublevel set for some \(c \in \R\).
    Every point \(\vek{y}\) on the boundary\footnote{The topological boundary of a set are all points such that every open neighbourhood of the point contains at least one point in the set and at least one point not in the set.} of \(\mathcal{L}_f\) or on the boundary of any connected component\footnote{The term \emph{connected component} refers to maximally connected subsets in the topological sense, see for example \cite[Section 4.5]{mendelson_introduction_1968}.} of \(\mathcal{L}_f\) satisfies \(f(\vek{y}) = c\).
\end{lemma}

\begin{proof}
    First, we show that \(f(\vek{y}) = c\) for all \(\vek{y} \in \partial \mathcal{L}_f\).
    Let \(\vek{z} \in \R^n\) be such that \(f(\vek{z}) > c\).
    By continuity, there exists a ball around \(\vek{z}\) with the property that the value of \(f\) is larger than \(c\) for every element of this ball.
    Therefore, \(\vek{z}\) lies in the interior of \(\R^n \setminus \mathcal{L}_f\) and not on the boundary of \(\mathcal{L}_f\).
    Similarly, every \(\vek{z} \in \R^n\) with \(f(\vek{z}) < c\) lies in the interior of \(\mathcal{L}_f\) and not on the boundary.
    This proves the first part.

    Now, consider any connected component \(\mathcal{C}\) of \(\mathcal{L}_f\).
    It suffices to show \(\partial \mathcal{C} \subseteq \partial \mathcal{L}_f\).
    For any point \(\vek{z} \notin \partial \mathcal{L}_f\) there exists a ball around \(\vek{z}\) such that this ball is either entirely contained in \(\mathcal{L}_f\) or entirely contained in \(\R^n \setminus \mathcal{L}_f\).
    In the latter case, the ball is also contained in \(\R^n \setminus \mathcal{C}\) and \(\vek{z}\) does not lie on the boundary \(\partial C\).
    In the former case, since every ball is connected, it is contained in one connected component of \(\mathcal{L}_f\).
    It is therefore either contained in \(\mathcal{C}\) or contained in \(\R^n \setminus \mathcal{C}\) and again \(\vek{z} \notin \partial C\).
\end{proof}

With these prerequisites, we can now prove that in particular local model minimizers close to \(\vek{x}_*\) are minimizers that will lead to a successful iteration.

\begin{lemma}\label{thm:close-minimizer}
    Let \(f \in C^p_{\Lip}\), \(\vek{x}_* \in M^q(f)\) and \(p > q-1\).
    There exists radii \(0 < r_x \leq r_y \leq r_q\) such that for any \(\vek{x}_k \in \ball(\vek{x}_*, r_x)\) there exists a local minimizer of \(m_k\) inside \(\ball(\vek{x}_*, r_y)\) and any approximate minimizer \(\vek{y}_k\) of \(m_k\) in the sense of \cref{eqn:approximate-minimizer} inside \(\ball(\vek{x}_*, r_y)\) has a success ratio \(\rho_k \geq \eta\).
\end{lemma}

\begin{proof}
    We first derive bounds for \(\vek{x}_k\) and \(\vek{y}_k\) such that \(\rho_k \geq \eta\), assuming that \(\vek{y}_k\) is an approximate minimizer of \(m_k\) satisfying \cref{eqn:approximate-minimizer}.
    It suffices to show that \(\abs{\rho_k - 1} \leq 1 - \eta\).
    Since
    \begin{equation}\label{eqn:rho-bound}
        \abs{\rho_k - 1} = \Abs{\frac{f(\vek{x}_k) - f(\vek{y}_k)}{t_k(\vek{x}_k) - t_k(\vek{y}_k)} - 1} = \frac{\abs{t_k(\vek{y}_k) - f(\vek{y}_k)}}{\abs{t_k(\vek{x}_k) - t_k(\vek{y}_k)}}
    \end{equation}
    we need to find an upper bound for the numerator and a lower bound for the denominator.
    The upper bound follows directly from \cref{eqn:taylor-error-bound}.

    Take the derivative of \cref{eqn:arp-model}, multiply with \(\vek{x}_k - \vek{y}_k\) and rearrange to obtain
    \begin{align*}
        \nabla t_k(\vek{y}_k)^T (\vek{x}_k - \vek{y}_k) &= \nabla m_k(\vek{y}_k)^T (\vek{x}_k - \vek{y}_k) + (p+1) \sigma_k \norm{\vek{y}_k - \vek{x}_k}^{p+1} \\
        &\geq -\norm{\nabla m_k(\vek{y}_k)} \norm{\vek{x}_k - \vek{y}_k} \geq -\theta \norm{\vek{x}_k - \vek{y}_k}^{p+1}.
    \end{align*}
    Here, the last inequality used that \(\vek{y}_k\) is an approximate minimizer in the sense of \cref{eqn:approximate-minimizer}.
    By \cref{thm:taylor-partial-uniform-convexity} there is a ball of radius \(r_q' > 0\) around \(\vek{x}_*\) such that \cref{eqn:taylor-uniform-convexity-x,eqn:taylor-uniform-convexity-y} hold inside that ball.
    Assume \(\norm{\vek{x}_k - \vek{y}_k} \leq \sqrt[p-q+1]{\frac{\mu_q'}{2 q \theta}}\) and \(\vek{x}_k, \vek{y}_k \in \ball(\vek{x}_*, r_q')\), then combining \cref{eqn:taylor-uniform-convexity-x} and hte previous inequality we get
    \begin{align*}
        t_k(\vek{x}_k) - t_k(\vek{y}_k) &\geq \nabla t_k(\vek{y}_k)^T (\vek{x}_k - \vek{y}_k) + \frac{\mu_q'}{q} \norm{\vek{x}_k - \vek{y}_k}^q \\
        &\geq \Paren{\frac{\mu_q'}{q} - \theta \norm{\vek{x}_k - \vek{y}_k}^{p-q+1}} \norm{\vek{x}_k - \vek{y}_k}^q \geq \frac{\mu_q'}{2q} \norm{\vek{x}_k - \vek{y}_k}^q.
    \end{align*}

    Therefore, the numerator in \cref{eqn:rho-bound} is bounded by an \(O(\norm{\vek{x}_k - \vek{y}_k}^{p+1})\) term and the denominator by an \(O(\norm{\vek{x}_k - \vek{y}_k}^q)\) term.
    This means
    \[
        \abs{\rho_k - 1} = \frac{\abs{t_k(\vek{y}_k) - f(\vek{y}_k)}}{\abs{t_k(\vek{x}_k) - t_k(\vek{y}_k)}} \leq \frac{L_p}{(p+1)!} \frac{2q}{\mu_q'} \norm{\vek{x}_k - \vek{y}_k}^{p-q+1} \leq 1 - \eta
    \]
    if additionally \(\norm{\vek{x}_k - \vek{y}_k} \leq \sqrt[p-q+1]{(1 - \eta) \frac{(p+1)!}{L_p} \frac{\mu_q'}{2q}}\).
    Overall, any approximate minimizer \(\vek{y}_k\) leads to a successful iteration whenever \(\vek{x}_k, \vek{y}_k \in \ball(\vek{x}_*, r_y)\) for
    \[
        r_y = \min\left\{r_q', \ \frac{1}{2} \sqrt[p-q+1]{\frac{\mu_q'}{2 q \theta}}, \ \frac{1}{2} \sqrt[p-q+1]{(1 - \eta) \frac{(p+1)!}{L_p} \frac{\mu_q'}{2q}}\right\}.
    \]

    What is left to show is that for \(\vek{x}_k \in \ball(\vek{x}_*, r_x)\) such a minimizer exists in the first place.
    When \(\vek{x}_k\) and \(\vek{y}\) are close enough to \(\vek{x}_*\) then, using \cref{eqn:taylor-uniform-convexity-y}, we can derive a lower bound for \(m_k\) as
    \begin{equation}\label{eqn:model-qth-order-lower-bound}
        m_k(\vek{y}) \geq t_k(\vek{y}) \geq t_k(\vek{x}_k) + \nabla t_k(\vek{x}_k)^T (\vek{y} - \vek{x}_k) + \frac{\mu_{q}'}{q} \norm{\vek{y} - \vek{x}_k}^q \eqqcolon \underline{m}_k(\vek{y}).
    \end{equation}
    Independently of whether the bound holds for any local minimizer of \(m_k\), observe that \(\underline{m}_k(\vek{x}_k) = m_k(\vek{x}_k)\), that \(\underline{m}_k\) is convex, and that the sublevel set
    \[
        \underline{\mathcal{L}}_k = \{\vek{y} \in \R^n \mid \underline{m}_k(\vek{y}) \leq m_k(\vek{x}_k) \}
    \]
    is contained in a closed ball of radius \(\sqrt[q-1]{\frac{q}{\mu_q'} \norm{\nabla t_k(\vek{x}_k)}}\) around \(\vek{x}_k\), because
    \begin{equation}\label{eqn:model-qth-order-lower-bound-large}
        \underline{m}_k(\vek{y}) \geq \underbrace{t_k(\vek{x}_k)}_{= m_k(\vek{x}_k)} + \underbrace{\Paren{- \norm{\nabla t_k(\vek{x}_k)} + \frac{\mu_q'}{q} \norm{\vek{y} - \vek{x}_k}^{q-1}}}_{> 0} \norm{\vek{y} - \vek{x}_k} > m_k(\vek{x}_k).
    \end{equation}
    for all \(\vek{y}\) outside of that ball.
    We need to make sure that \(\underline{\mathcal{L}}_k\) is contained in \(\ball(\vek{x}_*, r_q')\) such that \cref{eqn:model-qth-order-lower-bound} holds for points in the sublevel set.
    Let \(\vek{x}_k \in \ball(\vek{x}_*, r_q')\) and \(\vek{y} \in \underline{\mathcal{L}}_k\), then the distance of \(\vek{y}\) to \(\vek{x}_*\) can be bounded by the triangle inequality and
    \begin{align*}
        \norm{\vek{y} - \vek{x}_k} \leq \sqrt[q-1]{\frac{q}{\mu_q'} \norm{\nabla t_k(\vek{x}_k)}} \leq \sqrt[q-1]{\frac{q \nu}{\mu_q'} \norm{\vek{x}_k - \vek{x}_*}}
    \end{align*}
    since \(\norm{\nabla t_k(\vek{x}_k)} = \norm{\nabla f(\vek{x}_k)} \leq \nu \norm{\vek{x}_k - \vek{x}_*}\) by \cref{thm:uniform-convexity-properties}.
    For
    \begin{equation}\label{eqn:rx-definition}
        r_x = \min \left\{ r_q', \ \frac{\mu_q'}{q \nu} \Paren{\frac{r_y}{2}}^{q-1}, \ \frac{r_y}{2} \right\}
    \end{equation}
    and any expansion point \(\vek{x}_k \in \ball(\vek{x}_*, r_x)\) we know that any point \(\vek{y} \in \underline{\mathcal{L}}_k\) is contained in \(\ball(\vek{x}_*, r_y) \subseteq \ball(\vek{x}_*, r_q')\) which implies that \cref{eqn:model-qth-order-lower-bound} holds for all such points.
    Since \(\underline{m}_k\) is convex, \(\underline{\mathcal{L}}_k\) is a closed convex set and \(\vek{x}_k\) lies on its boundary.
    We always assume that the gradient of \(f\) at any iterate \(\vek{x}_k\) is nonzero, i.e.\ that \(\nabla \underline{m}_k(\vek{x}_k) = \nabla t_k(\vek{x}_k) = \nabla m_k(\vek{x}_k) = \nabla f(\vek{x}_k) \neq \vek{0}\).
    This implies that there are points \(\vek{y}\) in the interior of \(\underline{\mathcal{L}}_k\) where \(m_k(\vek{y}) < m_k(\vek{x}_k)\).
    By compactness, \(m_k\) attains its minimum on \(\underline{\mathcal{L}}_k\), and the minimizer is not on the boundary, since by \cref{thm:boundary-sublevel-set} any point \(\vek{y}\) on the boundary satisfies \(m_k(\vek{y}) \geq \underline{m}_k(\vek{y}) = m_k(\vek{x}_k)\).
    The point at which the minimum is attained must therefore lie in the interior of \(\underline{\mathcal{L}}_k\) and be a local minimizer of \(m_k\) which is contained in \(\ball(\vek{x}_*, r_y)\).
    This concludes the proof.
\end{proof}

The way of characterizing successful model minimizers in \cref{thm:close-minimizer} has a direct dependence on the function minimizer \(\vek{x}_*\), which the algorithm aims to find and is thus a priori unknown.
Alternatively, it is also possible to show that (asymptotically) it suffices to choose a local model minimizer from the correct connected component of the model sublevel set.
Clearly, by \cref{eqn:approximate-minimizer} the tentative iterate \(\vek{y}_k\) lies in the sublevel set \(\{\vek{y} \in \R^n \mid m_k(\vek{y}) \leq m_k(\vek{x}_k)\}\).
Since the local model does not have to be convex, there can be multiple connected components of this sublevel set.
Selecting a minimizer from the connected component that contains \(\vek{x}_k\) avoids spurious model minimizers far away from \(\vek{x}_*\) and ensures a successful step.

\begin{lemma}\label{thm:connected-component-minimizer}
    Let \(f \in C^p_{\Lip}\), \(\vek{x}_* \in M^q(f)\) and \(p > q-1\).
    Moreover, let \(\mathcal{C}_k\) be the connected component of the sublevel set \(\mathcal{L}_k \deq \{\vek{y} \in \R^n \mid m_k(\vek{y}) \leq m_k(\vek{x}_k)\}\) that contains \(\vek{x}_k\).
    There exists a radius \(r_c > 0\) such that for any \(\vek{x}_k \in \ball(\vek{x}_*, r_c)\) there exists a local minimizer of \(m_k\) inside \(\mathcal{C}_k\) and any approximate minimizer \(\vek{y}_k\) of \(m_k\) in the sense of \cref{eqn:approximate-minimizer} inside \(\mathcal{C}_k\) has a success ratio \(\rho_k \geq \eta\).
\end{lemma}

\begin{proof}
    The assumptions are the same as for \cref{thm:close-minimizer}, so there exists a radii \(r_x, r_y > 0\) such that all approximate minimizers \(\vek{y}_k\) inside \(\ball(\vek{x}_*, r_y)\) have a success ratio \(\rho_k \geq \eta\) whenever \(\vek{x}_k \in \ball(\vek{x}_*, r_x)\).
    Moreover, following the proof of that Lemma, we see that, on a ball of radius \(r_q' > 0\) around \(\vek{x}_*\), equation \cref{eqn:model-qth-order-lower-bound} holds and that \(\underline{m}_k\) is a convex lower bound for \(m_k\).
    By \cref{eqn:model-qth-order-lower-bound-large} this lower bound is strictly larger than \(m_k(\vek{x}_k)\) outside the ball \(\ball(\vek{x}_k, \sqrt[q-1]{\frac{q}{\mu_q'} \norm{\nabla t_k(\vek{x}_k)}})\).
    In other words,
    \[
        \mathcal{L}_k \cap \Paren{\ball(\vek{x}_*, r_q') \setminus \ball\Paren{\vek{x}_k, \sqrt[q-1]{\frac{q}{\mu_q'} \norm{\nabla t_k(\vek{x}_k)}}}} = \emptyset.
    \]

    Take \(r_c \in (0, r_x)\) with \(r_x\) as defined in \cref{eqn:rx-definition} and assume \(\vek{x}_k \in \ball(\vek{x}_*, r_c)\) from now on.
    The radius \(r_c\) is chosen small enough such that the previous equation implies \(\mathcal{L}_k \cap \partial \ball(\vek{x}_*, r_q') = \emptyset\).
    Thus, every connected component is either completely contained in the interior of \(\ball(\vek{x}_*, r_q')\) or completely contained in \(\R^n \setminus \ball(\vek{x}_*, r_q')\).
    Otherwise, the connected component could be partitioned into two sets, each the intersection of the component with an open set, which contradicts the connectedness of the connected component.
    Clearly, \(\mathcal{C}_k\) is contained in \(\ball(\vek{x}_*, r_q')\) and thus
    \begin{equation}\label{eqn:local-component-bounded-ry}
        \mathcal{C}_k \subseteq \ball\Paren{\vek{x}_k, \sqrt[q-1]{\frac{q}{\mu_q'} \norm{\nabla t_k(\vek{x}_k)}}} \subseteq \ball(\vek{x}_*, r_y).
    \end{equation}
    Every approximate minimizer \(\vek{y}_k \in \mathcal{C}_k\) must therefore satisfy \(\rho_k \geq \eta\).

    We now turn towards proving that there always exists a local minimizer of \(m_k\) in \(\mathcal{C}_k\).
    By \cref{eqn:local-component-bounded-ry} we know that \(\mathcal{C}_k\) is bounded and since every connected component of a closed set is closed \cite[Corollary 5.6, Chapter 4]{mendelson_introduction_1968}, \(\mathcal{C}_k\) must be compact.
    Therefore, the continuous function \(m_k\) attains its minimum on \(\mathcal{C}_k\) and by \cref{thm:boundary-sublevel-set} every point on the boundary of \(\mathcal{C}_k\) has the same function value as \(\vek{x}_k\).
    At the same time, the fact that \(\nabla m_k(\vek{x}_k) = \nabla f(\vek{x}_k) \neq \vek{0}\) implies the existence of a point with a lower function value in the interior of \(\mathcal{C}_k\).
    This means the minimum of \(m_k\) on \(\mathcal{C}_k\) is attained in the interior of \(\mathcal{C}_k\) and so this point is a local minimizer of \(m_k\).
\end{proof}

\begin{theorem}\label{thm:convergence-pq-rate-good-min}
    Let \(f \in C^p_{\Lip}\), \(\vek{x}_* \in M^q(f)\) and \(p > q-1\).
    Assume that \(\vek{y}_k\) is always chosen such that \(\vek{y}_k \in \ball(\vek{x}_*, r_y)\), with \(r_y\) as defined in \cref{thm:close-minimizer}, or always chosen in a way that \(\vek{y}_k \in \mathcal{C}_k\), with \(\mathcal{C}_k\) as defined in \cref{thm:connected-component-minimizer}.
    If additionally, \(\vek{x}_0\) is close enough to \(\vek{x}_*\), then all iterations of \cref{alg:arp} are successful and
    \begin{align*}
        \norm{\nabla f(\vek{x}_k)} &\to 0 && \text{at a Q-} \tfrac{p}{q-1} \text{th-order rate} \\
        f(\vek{x}_k) &\to f(\vek{x}_*) && \text{at a Q-} \tfrac{p}{q-1} \text{th-order rate} \\
        \vek{x}_k &\to \vek{x}_* && \text{at an R-} \tfrac{p}{q-1} \text{th-order rate}
    \end{align*}
    as \(k \to \infty\).
\end{theorem}

\begin{proof}
    Note that in the proof of \cref{thm:connected-component-minimizer}, the radius \(r_c > 0\) is chosen small enough such that for any \(\vek{x}_k \in \ball(\vek{x}_*, r_c)\) the component satisfies \(\mathcal{C}_k \subseteq \ball(\vek{x}_*, r_y)\).
    Therefore, for \(\vek{x}_k\) close enough to \(\vek{x}_*\) the assumption that \(\vek{y}_k \in \ball(\vek{x}_*, r_y)\) is the weaker one, and it suffices to prove the result for this assumption.

    Let \(r_x\) and \(r_y\) be as defined in \cref{thm:close-minimizer}.
    We will first need to show inductively that, for \(\vek{x}_0\) close enough to \(\vek{x}_*\), all \(\vek{x}_k \in \ball(\vek{x}_*, r_x)\).\footnote{For the \(\vek{y}_k \in \mathcal{C}_k\) assumption based on \cref{thm:connected-component-minimizer} we would of course need to prove that all iterates stay in \(\ball(\vek{x}_*, r_c)\). The argument is the same after replacing \(r_x\) with \(r_c\).}
    Choose \(\vek{x}_0 \in \ball(\vek{x}_*, r_x)\) close enough to \(\vek{x}_*\) such that
    \[
        \norm{\nabla f(\vek{x}_0)} \leq \min \left\{ \Paren{\frac{1}{R_{\max}} \Paren{\frac{\mu_q}{q}}^{\frac{p}{q-1}}}^{\frac{q-1}{p-q+1}}, \ \mu_q r_x^{q-1} \right\}
    \]
    For the induction step assume \(\vek{x}_k \in \ball(\vek{x}_*, r_x)\) is such that \(\norm{\nabla f(\vek{x}_k)} \leq \norm{\nabla f(\vek{x}_0)}\).
    By \cref{thm:close-minimizer} there exists an approximate model minimizer \(\vek{y}_k \in \ball(\vek{x}_*, r_y)\) and by assumption such a point is chosen.
    This means that the iteration is successful and \(\vek{x}_k, \vek{x}_{k+1} = \vek{y}_k \in \ball(\vek{x}_*, r_q)\), i.e., both points lie in the region where \(f\) is uniformly convex.
    Apply \cref{thm:gradient-convergence-pq-rate} to obtain
    \[
        \norm{\nabla f(\vek{x}_{k+1})} \leq R_{\max} \Paren{\frac{q}{\mu_q}}^{\frac{p}{q-1}} \norm{\nabla f(\vek{x}_k)}^{\frac{p-q+1}{q-1}} \norm{\nabla f(\vek{x}_k)},
    \]
    which implies \(\norm{\nabla f(\vek{x}_{k+1})} \leq \norm{\nabla f(\vek{x}_k)} \leq \norm{\nabla f(\vek{x}_0)}\).
    Therefore, by \cref{eqn:uniform-convexity-step-bounded-by-gradient}
    \[
        \norm{\vek{x}_{k+1} - \vek{x}_*} \leq \sqrt[q-1]{\frac{\norm{\nabla f(\vek{x}_{k+1})}}{\mu_q}} \leq \sqrt[q-1]{\frac{\norm{\nabla f(\vek{x}_0)}}{\mu_q}} \leq r_x.
    \]
    We conclude that \(\vek{x}_{k+1} \in \ball(\vek{x}_*, r_x)\) and satisfies \(\norm{\nabla f(\vek{x}_{k+1})} \leq \norm{\nabla f(\vek{x}_0)}\), thus completing the induction.
    As mentioned before, all iterations are successful and the rates follow directly from \cref{thm:convergence-pq-rate-successful}.
\end{proof}

To conclude this section, we comment on different aspects of the results above.
First, \cref{thm:convergence-pq-rate-good-min} implies that for fully adaptive methods there is a way to ensure that asymptotically all iterations are successful.
However, the way to do so is not, as one might expect, to choose the global model minimizer, as \cref{exm:sigma-oscillations} shows.
Instead, when faced with ambiguity one should select a local minimizer of the model that lies in the vicinity of \(\vek{x}_*\).
While neither \cref{thm:close-minimizer} nor \cref{thm:connected-component-minimizer} give characterizations of the right model minimizer that lead to practical algorithms for finding such minimizers, \cref{thm:connected-component-minimizer} suggests that any monotonically decreasing local optimization algorithm applied to \(m_k\) and starting at \(\vek{x}_k\) is unlikely compute a wrong minimizer.
To do so such a local method would need to compute a large enough step in the right direction to jump over a hill in \(m_k\) in order to move from \(\mathcal{C}_k\) to a different connected component of the sublevel set.

Second, note that neither in the statement nor in the proof of \cref{thm:close-minimizer} there is an assumption or dependence on \(\sigma_k\), other that \(\sigma_k \geq 0\).
The existence of this local model minimizer close to \(\vek{x}_*\) is even true for \(\sigma_k = 0\).
Thus, when selecting \(\vek{y}_k\) in the appropriate way, it is possible to drop the regularization parameter to zero once \(\vek{x}_k\) is close enough to \(\vek{x}_*\) and still achieve the same convergence guarantees.
Unfortunately, without regularization the subproblem only reduces to solving a single linear system in the second-order case.
For third- and higher-order algorithms even minimizing the unregularized Taylor expansion is often done with general local optimization methods and gets more complicated as \(p\) increases.

Lastly, the claims and proofs of this section can easily be extended from functions with Lipschitz continuous \(p\)th derivatives to ones with H\"older continuous \(p\)th derivatives, i.e. those that satisfy
\[
    \norm{\nabla^p f(\vek{x}) - \nabla^p f(\vek{y})} \leq H_p \norm{\vek{x} - \vek{y}}^{\beta_p}
\]
for some constants \(H_p > 0\), \(\beta_p \in [0, 1]\) for all \(\vek{x}, \vek{y} \in \R^n\).
Let AR\(p, r\) be the variant of \cref{alg:arp} where the exponent of the regularization term is \(r\) and the exponent in the conditions for an approximate local minimizer \cref{eqn:approximate-minimizer} is \(r-1\).
Applying AR\(p, p+\beta_p\) to functions that satisfy the conditions above gives \(\frac{p + \beta_p - 1}{q - 1}\)-order convergence (compared to \(\frac{p}{q-1}\) in \cref{thm:convergence-pq-rate-successful}) whenever \(p + \beta_p > q\).
\Cref{thm:convergence-pq-rate,thm:convergence-pq-rate-good-min} change correspondingly.

From one perspective, the analysis under the assumption of H\"older continuity bridges the gaps and clarifies how smoothness (\(p + \beta_p\)) and convexity (\(q\)) interact.
With the restriction to Lipschitz continuity the convergence rates seem to be improving in discrete steps as \(p\) increases, while the result above makes clear that actually these rates are part of a continuum, where the rates continuously improve as \(p + \beta_p\) increases.
Unfortunately, the \(\frac{p + \beta_p - 1}{q - 1}\)-order rate is only possibly with a priori knowledge of the \(\beta_p\) parameter as the value of \(\beta_p\) is used within the algorithm.\footnote{
    An analysis of a variant of the adaptive regularization method when the regularization power \(r\) and the smoothness \(p + \beta_p\) are not equal is provided by Cartis, Gould and Toint \cite{cartis_universal_2019}.
    Their results show that the bound on the number of iterations does not degrade when \(r \geq p + \beta_p\) and varies continuously with \(p + \beta_p\).
    With this reasoning the authors suggest that \(r = p+1\) gives a universal algorithm that adapts to the parameter \(\beta_p\) of the function class at hand.
    Whether \((p+1)\)st-order regularization also leads to a universal algorithm in the context of local convergence, is an open question and outside the scope of this paper.
}

From a second perspective, the Lipschitz continuous case is the most important case and the order of convergence does indeed improve in discrete steps as the algorithm gets access to more and more derivatives.
When considering smooth functions \(f \in C^{\infty}\) and applying AR\(p\) for different \(p\), then these functions have Lipschitz continuous \(p\)th derivatives on any bounded domain.
With the additional assumption that the iterates stay bounded, we obtain rates that depend not on the smoothness itself, but on the number of derivatives the algorithm has access to, which can only be increased in discrete steps.
This is why in this paper we focus on the Lipschitz continuity case and only mention the more general H\"older continuity case here.

\section{Sharpness of local convergence results}\label{sec:sharpness}

Considering the \(\frac{p}{q-1}\)th-order convergence rates for successful iterations in \cref{thm:convergence-pq-rate-successful} and the degraded rates in \cref{thm:convergence-pq-rate}, it is natural to ask whether these rates are sharp and in which cases one might need to expect such rates.
The example in \cref{sec:sharp-pq-convergence} shows that the order of convergence proved in \cref{thm:convergence-pq-rate-successful} is indeed sharp and the main result in \cref{sec:sharp-oscillations} proves that in some cases oscillations in \(\sigma_k\) are unavoidable when \(\vek{y}_k\) is the global model minimizer and \(\sigma_k\) is decreased during successful iterations.

\subsection{\texorpdfstring{\(\frac{p}{q-1}\)th-order convergence}{p/(q-1)th-order convergence}}\label{sec:sharp-pq-convergence}

\begin{example}\label{exm:sharp-pq-convergence}
    Consider the one-dimensional function \(f(x) = \frac{1}{q} x^q + \frac{1}{p+1} x^{p+1}\) for even \(q \in \N\) and some \(p \in \N\) with \(p > q-1\).
    The function has a local minimizer at \(x_* = 0\), its \(p\)th derivative is Lipschitz continuous with constant \(L_p = p!\), and it is locally uniformly convex of order \(q\).
    We choose \(\theta = 3\) and the other parameters arbitrarily.
    Since \(p \geq q\), the \(p\)th-order Taylor expansion of \(x^q\) is exact and the expansion of \(x^{p+1}\) is only missing the \((p+1)\)st-order term.
    Therefore,
    \begin{subequations}
        \begin{align}
            t_k(y) &= f(y) - \frac{1}{p+1} (y - x_k)^{p+1} = \frac{1}{q} y^q + \frac{1}{p+1} y^{p+1} - \frac{1}{p+1} (y - x_k)^{p+1} \\
            t_k'(y) &= y^{q-1} + y^p - (y - x_k)^p \\
            m_k'(y) &= y^{q-1} + y^p - (y - x_k)^p + (p+1) \sigma_k \abs{y - x_k}^{p-1} (y - x_k). \label{eqn:pq-example-model-derivative}
        \end{align}
    \end{subequations}
    At each iteration, the AR\(p\) algorithm computes an approximate minimizer \(y_k\) of \(m_k\) satisfying \cref{eqn:approximate-minimizer}.
    We will show that
    \[
        y_k = - \abs{x_k}^{\frac{p}{q-1}} \frac{x_k}{\abs{x_k}}
    \]
    is an approximate minimizer as long as \(\abs{x_k}\) is small enough.
    Notice that the sign of \(y_k\) is always the opposite of the sign of \(x_k\), which implies \(\abs{y_k - x_k} = \abs{y_k} + \abs{x_k}\).
    We can deduce that
    \[
        \abs{y_k^{q-1} + y_k^p} = \Abs{-\abs{x_k}^p \frac{x_k}{\abs{x_k}} + y_k^p} \leq \abs{x_k}^p + \abs{y_k}^p \leq \Paren{\abs{x_k} + \abs{y_k}}^p = \abs{y_k - x_k}^p.
    \]
    The other two terms in \cref{eqn:pq-example-model-derivative} are bounded by \(\abs{y_k - x_k}^p\) and \((p+1)\sigma_k \abs{y_k - x_k}^p\) respectively. Overall,
    \[
        \abs{m_k'(y_k)} \leq (2 + (p+1) \sigma_k) \abs{y_k - x_k}^p \leq \theta \abs{y_k - x_k}^p
    \]
    for \(\sigma_k \leq \frac{1}{p+1}\).
    This makes the tentative iterate \(y_k\) an approximate minimizer since \(m_k(y_k) < m_k(x_k)\) as long as \(\abs{x_k}\) is small enough.
    Moreover, we can make \(\abs{x_k}\) small enough such that \cref{thm:close-minimizer} implies that any approximate minimizer in a prescribed neighbourhood of \(x_* = 0\) is successful and that \(y_k\), satisfying \(\abs{y_k} = \abs{x_k}^{\frac{p}{q-1}}\), lies in this neighbourhood.
    For \(\sigma_0 < \frac{1}{p+1}\) and \(\abs{x_0}\) small enough, all iterations are thus successful and
    \[
        \abs{x_k} = \abs{x_0}^{\Paren{\frac{p}{q-1}}^k}.
    \]
    Therefore, \(x_k\) converges to \(x_*\) at a \(\frac{p}{q-1}\)th-order rate, which in turn implies that the order of convergence of \(f'(x_k) \to 0\) and \(f(x_k) \to f(x_*)\) is not higher by \cref{thm:uniform-convexity-properties} (but also not lower by \cref{thm:convergence-pq-rate-successful}).
\end{example}

This example illustrates that the convergence rates derived in \cref{thm:convergence-pq-rate-successful,thm:convergence-pq-rate-good-min} are not an artefact of the analysis, but describe worst-case rates that can be attained.

\subsection{Unsuccessful iterations}\label{sec:sharp-oscillations}

\Cref{exm:sigma-oscillations} already shows that there are cases where half of the iterations are unsuccessful even asymptotically, when \(\vek{y}_k\) is the global model minimizer at each iteration and \(\sigma_k\) is decreased during successful iterations.
Continuing with this global minimizer assumption, we want to make more precise in which cases we need to expect unsuccessful iterations, show that they occur for a wide variety of problems, and that the ratio of unsuccessful iterations to successful iterations is equal to \(\alpha \geq 0\) where \(\alpha\) describes the relationship between the increase and decrease factors for \(\sigma_k\) through \(\gamma_1 = \gamma_2^{-\alpha}\).
This then shows that the degradation from \(\frac{p}{q-1}\)th-order convergence rates for successful iterations to \(\sqrt[1+\alpha]{\frac{p}{q-1}}\)th-order rates in \cref{thm:convergence-pq-rate} are again sharp.

The key assumption in the following theorem is that \(\vek{x}_*\) is a global minimizer of \(f\) but not of \(t_{\vek{x}_*}\), the Taylor expansion at \(\vek{x}_*\).
The second part holds true, for example, when \(p\) is odd and the highest-order term is non-zero, because then the Taylor expansion is not lower bounded and has no global minimizer.
When these two assumptions are satisfied, the smallest \(\sigma_k\) that guarantees a successful iteration and the largest \(\sigma_k\) that guarantees an unsuccessful iteration both converge to some value \(\sigma_*\) (\cref{thm:oscillations-sigma-bounds}).
Therefore, the regularization parameters \(\sigma_k\) will eventually oscillate around \(\sigma_*\) (\cref{thm:oscillations}).

To prove this result we need some bound on the difference between the model functions when expanded at \(\vek{x}_k\) and \(\vek{x}_*\) to then exploit the fact that \(\vek{x}_*\) is a local but not a global minimizer of the Taylor expansion at \(\vek{x}_*\).

\begin{lemma}\label{thm:model-difference-bound}
    Let \(f \in C^p_{\Lip}\).
    For any \(r > 0\) and \(\vek{x} \in \R^n\) there exists a constant \(C_{\vek{x}, r} > 0\) such that
    \begin{subequations}\label{eqn:model-difference-bound}
        \begin{align}
            \abs{m_{\vek{x}, \sigma}(\vek{y}) - m_{\vek{x}', \sigma}(\vek{y})} &\leq (C_{\vek{x}, r} + 2^p (p+1) \sigma) \norm{\vek{x} - \vek{x}'} \norm{\vek{y} - \vek{x}}^p \text{ and} \\
            \abs{m_{\vek{x}, \sigma}(\vek{y}) - m_{\vek{x}', \sigma}(\vek{y})} &\leq (C_{\vek{x}, r} + 2^p (p+1) \sigma) \norm{\vek{x} - \vek{x}'} \norm{\vek{y} - \vek{x}'}^p
        \end{align}
    \end{subequations}
    hold for all \(\vek{x}' \in \ball(\vek{x}, r)\) and \(\vek{y} \notin \ball(\vek{x}, 2r)\).
\end{lemma}

\begin{proof}
    To prove this bound on the difference between models with different expansion points, we need to first establish some simple bounds between differences of rank-one tensors and differences between \((p+1)\)st-order regularization terms.
    For any \(\vek{x}, \vek{x}', \vek{y} \in \R^n\) and \(j \in \N_{\geq 0}\) we have
    \begin{equation}\label{eqn:rank-one-tensor-difference}
        \norm{\otimes^j (\vek{y} - \vek{x})^T - \otimes^j (\vek{y} - \vek{x}')^T} \leq j \norm{\vek{x} - \vek{x}'} \max(\norm{\vek{y} - \vek{x}}, \norm{\vek{y} - \vek{x}'})^{j-1}.
    \end{equation}
    This follows from the following telescoping sum
    \[
        \phantom{{}={}} \otimes^j (\vek{y} - \vek{x})^T - \otimes^j (\vek{y} - \vek{x}')^T = \sum_{i=0}^{j-1} \otimes^{j-i-1} (\vek{y} - \vek{x})^T \otimes (\vek{x}' - \vek{x})^T \otimes^i (\vek{y} - \vek{x}')^T
    \]
    by taking norms and bounding each term individually.
    Note that when \(j=0\) the ``tensors'' involved are outer products of zero vectors, i.e. the scalar \(1\), so the difference between them is zero and so is the right-hand side of \cref{eqn:rank-one-tensor-difference}.
    In the same way,
    \[
        \Abs{\norm{\vek{y} - \vek{x}}^{p+1} - \norm{\vek{y} - \vek{x}'}^{p+1}} \leq (p+1) \norm{\vek{x} - \vek{x}'} \max(\norm{\vek{y} - \vek{x}}, \norm{\vek{y} - \vek{x}'})^p
    \]
    can be derived from a telescoping sum using the reverse triangle inequality \(\abs{\norm{\vek{y} - \vek{x}} - \norm{\vek{y} - \vek{x}'}} \leq \norm{\vek{x} - \vek{x}'}\).

    Let us now turn to proving the original statement, by considering the difference between the Taylor expansions.
    \[
        t_{\vek{x}}(\vek{y}) - t_{\vek{x}'}(\vek{y}) = \sum_{j=0}^p \frac{1}{j!} \Paren{\Innerprod{\nabla^j f(\vek{x})}{\otimes^j (\vek{y} - \vek{x})^T}_F - \Innerprod{\nabla^j f(\vek{x}')}{\otimes^j (\vek{y} - \vek{x}')^T}_F}
    \]
    Ignoring the \(\frac{1}{j!}\) coefficient, each summand on the right-hand side can be written as
    \begin{equation}\label{eqn:taylor-expansion-difference-term-bound}
        \Innerprod{\nabla^j f(\vek{x})}{\otimes^j (\vek{y} - \vek{x})^T - \otimes^j (\vek{y} - \vek{x}')^T}_F + \Innerprod{\nabla^j f(\vek{x}) - \nabla^j f(\vek{x}')}{\otimes^j (\vek{y} - \vek{x}')^T}_F.
    \end{equation}
    For \(0 \leq j \leq p\) let \(M_{\vek{x}, j} = \norm{\nabla^j f(\vek{x})}\) and let \(L_{\vek{x}, r, j}\) be the local Lipschitz constant of \(\nabla^j f\) on \(\ball(\vek{x}, r)\).
    The local Lipschitz constants exist and are finite because \(\nabla^{j+1} f\) is continuous for \(0 \leq j \leq p-1\) and because \(\nabla^p f\) is Lipschitz continuous by assumption.
    Using these constants, the term in \cref{eqn:taylor-expansion-difference-term-bound} is bounded by
    \[
        M_{\vek{x}, j} j \norm{\vek{x} - \vek{x}'} \max(\norm{\vek{y} - \vek{x}}, \norm{\vek{y} - \vek{x}'})^{j-1} + L_{\vek{x}, r, j} \norm{\vek{x} - \vek{x}'} \max(\norm{\vek{y} - \vek{x}}, \norm{\vek{y} - \vek{x}'})^j.
    \]
    Thus, overall we get
    \begin{subequations}\label{eqn:model-difference-bound-all-y}
        \begin{align}
            &\phantom{{}={}} \abs{m_{\vek{x}, \sigma}(\vek{y}) - m_{\vek{x}', \sigma}(\vek{y})} \\
            &\leq \abs{t_{\vek{x}}(\vek{y}) - t_{\vek{x}'}(\vek{y})} + (p+1) \sigma \norm{\vek{x} - \vek{x}'} \max(\norm{\vek{y} - \vek{x}}, \norm{\vek{y} - \vek{x}'})^p \\
            &\leq \sum_{j=0}^p \frac{L_{\vek{x}, r, j} + M_{\vek{x}, j+1}}{j!} \norm{\vek{x} - \vek{x}'} \max(\norm{\vek{y}-\vek{x}}, \norm{\vek{y} - \vek{x}'})^j \\
            &\phantom{{}={}} + (p+1) \sigma \norm{\vek{x} - \vek{x}'} \max(\norm{\vek{y} - \vek{x}}, \norm{\vek{y} - \vek{x}'})^p \notag{}
        \end{align}
    \end{subequations}
    where \(M_{\vek{x}, p+1}\) is defined to be zero.
    This result holds for all \(\vek{y} \in \R^n\).

    The simplified bounds in the statement of this lemma can be obtained using \(\vek{y} \notin \ball(\vek{x}, 2r)\), since then
    \begin{align*}
        \norm{\vek{y} - \vek{x}'} &\leq \norm{\vek{y} - \vek{x}} + \norm{\vek{x} - \vek{x}'} \leq \norm{\vek{y} - \vek{x}} + r < 2 \norm{\vek{y} - \vek{x}} \text{ and} \\
        \norm{\vek{y} - \vek{x}} &< \norm{\vek{y} - \vek{x}} + \norm{\vek{y} - \vek{x}} - 2 \norm{\vek{x} - \vek{x}'} = 2 \Paren{\norm{\vek{y} - \vek{x}} - \norm{\vek{x} - \vek{x}'}} \leq 2 \norm{\vek{y} - \vek{x}'}. \notag{}
    \end{align*}
    Using \(\max(\norm{\vek{y} - \vek{x}}, \norm{\vek{y} - \vek{x}'}) > 2r\) and the first of these inequalities gives
    \begin{equation}\label{eqn:y-distance-bound-simplification}
        \max(\norm{\vek{y} - \vek{x}}, \norm{\vek{y} - \vek{x}'})^j \leq (2r)^{j-p} \max(\norm{\vek{y} - \vek{x}}, \norm{\vek{y} - \vek{x}'})^p \leq 2^j r^{j-p} \norm{\vek{y} - \vek{x}}^p
    \end{equation}
    for all \(0 \leq j \leq p\), and analogously the second inequality implies that the left-hand side is upper-bounded by \(2^j r^{j-p} \norm{\vek{y} - \vek{x}'}^p\) as well.

    Lastly, combine \cref{eqn:model-difference-bound-all-y,eqn:y-distance-bound-simplification} and choose
    \[
        C_{\vek{x}, r} = \sum_{j=0}^p \frac{2^j}{r^{p-j}} \frac{L_{\vek{x}, r, j} + M_{\vek{x}, j+1}}{j!}
    \]
    to obtain \cref{eqn:model-difference-bound}.
    This proves the claim.
\end{proof}

\begin{theorem}\label{thm:oscillations-sigma-bounds}
    Let \(f \in C^p_{\Lip}\), \(\vek{x}_* \in M^q(f)\) and \(p > q-1\).
    Assume that \(\vek{x}_*\) is a global minimizer of \(f\) but not a global minimizer of \(t_{\vek{x}_*}^p\).
    There exists a value \(\sigma_* > 0\) (independent of \(\eta\), \(\gamma_1\), and \(\gamma_2\)) such that for any \(\e > 0\) there exists a radius \(r_{\e} > 0\) with the property that for any \(\vek{x}_k \in \ball(\vek{x}_*, r_{\e})\) the iteration will be unsuccessful if \(\sigma_k < \sigma_* - \e\) and successful if \(\sigma_k > \sigma_* + \e\), given that \(\vek{y}_k\) is the global minimizer of \(m_k\).
\end{theorem}

\begin{proof}
    We will show that the claim holds for the following threshold value of the regularization parameter:
    \[
        \sigma_* = \inf \Set{\sigma > 0 | \min_{\vek{y} \in \R^n} m_{\vek{x}_*, \sigma}(\vek{y}) = f(\vek{x}_*)}.
    \]
    Since \(\vek{x}_*\) is not a global minimizer of \(t_{\vek{x}_*}^p\) there exists a point \(\vek{y}\) such that \(t_{\vek{x}_*}^p(\vek{y}) < t_{\vek{x}_*}^p (\vek{x}_*) = f(\vek{x}_*)\).
    If \(\sigma\) is small enough, the regularization term \(\sigma \norm{\vek{y} - \vek{x}_*}^p\) will be smaller than the difference between the two values of the Taylor expansion, therefore \(\sigma_* > 0\).
    \Cref{eqn:taylor-error-bound} implies that with \(\sigma \geq L_p / (p+1)!\) the model \(m_{\vek{x}_*, \sigma}\) overestimates the objective function everywhere.
    Combining this fact with \(m_{\vek{x}_*, \sigma}(\vek{x}_*) = f(\vek{x}_*)\) and the assumption that \(\vek{x}_*\) is a global minimizer of \(f\) we get \(\sigma_* < \infty\).

    Let \(\e > 0\) be given.
    We consider the case \(\sigma_k < \sigma_* - \e\) first.
    By definition of \(\sigma_*\) we know that there exists a point \(\vek{y}_{\e}\) with the property that \(m_{\vek{x}_*, \sigma_* - \e}(\vek{y}_{\e}) < f(\vek{x}_*)\).
    Let the difference between these two values be \(c\) and assume \(\norm{\vek{x}_k - \vek{x}_*} \leq r < \norm{\vek{y}_{\e} - \vek{x}_*}/2\), then the global minimizer \(\vek{y}_k\) of \(m_k\) satisfies
    \begin{align*}
        m_k(\vek{y}_k) &\leq m_k(\vek{y}_{\e}) = m_{\vek{x}_k, \sigma_k}(\vek{y}_{\e}) < m_{\vek{x}_k, \sigma_* - \e}(\vek{y}_{\e}) \\
        &\leq m_{\vek{x}_*, \sigma_* - \e}(\vek{y}_{\e}) + (C_{\vek{x}_*, r} + 2^p (p+1) (\sigma_* - \e)) \norm{\vek{x}_k - \vek{x}_*}\norm{\vek{y}_{\e} - \vek{x}_*}^p \\
        &\leq m_{\vek{x}_*, \sigma_* - \e}(\vek{y}_{\e}) + \tfrac{c}{2} \leq f(\vek{x}_*) - \tfrac{c}{2}
    \end{align*}
    if \(\norm{\vek{x}_k - \vek{x}_*}\) is small enough using \cref{thm:model-difference-bound}.
    This implies
    \[
        \rho_k = \frac{f(\vek{x}_k) - f(\vek{y}_k)}{t_k(\vek{x}_k) - t_k(\vek{y}_k)} \leq \frac{f(\vek{x}_k) - f(\vek{x}_*)}{f(\vek{x}_k) - m_k(\vek{y}_k)} \leq \frac{f(\vek{x}_k) - f(\vek{x}_*)}{f(\vek{x}_k) - f(\vek{x}_*) + \frac{c}{2}} < \eta
    \]
    whenever \(f(\vek{x}_k) - f(\vek{x}_*) < \frac{\eta}{1 - \eta} \frac{c}{2}\).
    Overall, as long as \(\norm{\vek{x}_k - \vek{x}_*}\) is small enough (which also means that \(f(\vek{x}_k) - f(\vek{x}_*)\) becomes arbitrarily small) any iteration with \(\sigma_k < \sigma_* - \e\) will be unsuccessful.

    Next, consider the case \(\sigma_k > \sigma_* + \e\).
    Choose \(r_x\) and \(r_y\) such that they satisfy the properties in \cref{thm:close-minimizer} and let \(r_x < r_y/2\).
    Apply \cref{thm:model-difference-bound} to \(\vek{x}_k \in \ball(\vek{x}_*, r_x)\) and \(\vek{y} \notin \ball(\vek{x}_*, r_y)\) to see that
    \begin{align*}
        &\phantom{{}={}} m_k(\vek{y}) = m_{\vek{x}_k, \sigma_k}(\vek{y}) > m_{\vek{x}_k, \sigma_*}(\vek{y}) + \e \norm{\vek{y} - \vek{x}_k}^{p+1} \\
        &\geq m_{\vek{x}_*, \sigma_*}(\vek{y}) + \Paren{\e \norm{\vek{y} - \vek{x}_k} - \Paren{C_{\vek{x}_*, r_x} + 2^p (p+1) \sigma_*} \norm{\vek{x}_k - \vek{x}_*}} \norm{\vek{y} - \vek{x}_k}^p.
    \end{align*}
    Note that \(\norm{\vek{y} - \vek{x}_k} \geq \norm{\vek{y} - \vek{x}_*} - \norm{\vek{x}_k - \vek{x}_*} > r_y/2\).
    Assuming
    \[
        \norm{\vek{x}_k - \vek{x}_*} \leq \frac{\e r_y}{4\Paren{C_{\vek{x}_*, r_x} + 2^p (p+1) \sigma_*}} \text{ and } f(\vek{x}_k) \leq f(\vek{x}_*) + \frac{\e}{2} \Paren{\frac{r_y}{2}}^{p+1}
    \]
    we can deduce that
    \begin{align*}
        m_k(\vek{y}) &> m_{\vek{x}_*, \sigma_*}(\vek{y}) + \frac{\e}{2} \frac{r_y}{2} \norm{\vek{y} - \vek{x}_k}^p \geq f(\vek{x}_*) + \frac{\e}{2} \Paren{\frac{r_y}{2}}^{p+1} \geq f(\vek{x}_k) = m_k(\vek{x}_k).
    \end{align*}
    Therefore, if \(\vek{x}_k\) is close enough to \(\vek{x}_*\) then any point outside \(\ball(\vek{x}_*, r_y)\) is not a global minimizer of the local model \(m_k\).
    Since we assumed \(\vek{y}_k\) is the global minimizer of \(m_k\) it must be contained in \(\ball(\vek{x}_*, r_y)\) and the iteration is successful by \cref{thm:close-minimizer}.
\end{proof}

\begin{theorem}\label{thm:oscillations}
    Let \(f \in C^p_{\Lip}\), \(\vek{x}_* \in M^q(f)\) and \(p > q-1\).
    Assume that \(\vek{x}_*\) is a global minimizer of \(f\) but not a global minimizer of \(t_{\vek{x}_*}^p\), that \(\gamma_1 = \gamma_2^{-\alpha}\) for some \(\alpha \in \Q\) and that \(\sigma_0 \neq \sigma_*\) for the value of \(\sigma_*\) from \cref{thm:oscillations-sigma-bounds}.
    For \(\vek{x}_0\) close enough to \(\vek{x}_*\) the values of \(\sigma_k\) will eventually cycle through finitely many values and the ratio of unsuccessful to successful iterations is \(\alpha\) in each cycle, if \(\vek{y}_k\) is the global minimizer of \(m_k\) at each iteration.
\end{theorem}

\begin{proof}
    By \cref{thm:oscillations-sigma-bounds} there exists the value \(\sigma_*\) for the given \(f\) and global minimizer \(\vek{x}_*\).
    Let \(\alpha = \frac{a}{b}\) for some \(a, b \in \N_{\geq 0}\) coprime if \(\alpha > 0\) and \(a = 0\) and \(b = 1\) if \(\alpha = 0\).
    Select
    \[
        \e < \min \Set{ \abs{\gamma_2^{c/b} \sigma_0 - \sigma_*} | c \in \Z }
    \]
    and let \(r_{\e}\) be chosen such that the conclusions of \cref{thm:oscillations-sigma-bounds} holds.
    By assuming \(\vek{x}_*\) to be a global minimizer, a small adaptation of \cref{thm:non-escape} implies that all iterates satisfy \(\vek{x}_k \in \ball(\vek{x}_*, r_{\e})\) if \(\vek{x}_0\) is close enough to \(\vek{x}_*\).
    Clearly, for any integer \(c\) we have \(\gamma_2^{c/b} \sigma_0 \notin [\sigma_* - \e, \sigma_* + \e]\).
    By the mechanism by which \(\sigma_k\) is updated in \cref{alg:arp}, we have that for all iterations \(\sigma_k = \gamma_2^{c_k/b} \sigma_0\) for some integer \(c_k \in \Z\).
    Let \(c_* \in \Z\) be such that \(\gamma_2^{c_*/b} \sigma_0\) is the smallest value larger than \(\sigma_*\).
    Whenever \(c_k \geq c_*\) the iteration is successful by \cref{thm:oscillations-sigma-bounds} and \(c_{k+1} = c_k - a\), and by the same logic whenever \(c_k < c_*\) the iteration is unsuccessful and \(c_{k+1} = c_k + b\).

    If \(a = 0\) then the regularization parameter will increase until \(c_k \geq c_*\) and then stay constant, since from that point onwards all iterations are successful.
    The claim holds as \(\sigma_k\) eventually cycles through one value and the ratio of unsuccessful to successful iterations is \(0 = \frac{a}{b}\).

    Assume \(a \neq 0\) from now on.
    Clearly, after each \(c_k \geq c_*\) there will eventually be a \(c_k < c_*\) and vice versa.
    Moreover, \(c_k\) will at some point enter the range \([c_* - a, c_* + b)\) and will never leave it again, so it suffices to show that once it does, the values cycle through finitely many values as claimed.

    Assume for the moment that \(c_0 = c_*\).
    Since there are only \(a + b\) distinct values that \(c_k\) can take, the values must cycle.
    Consider the remainder of \(c_k - c_*\) modulo \(b\), then it does not change in unsuccessful iterations \(c_k < c_*\), but it does change by \(-a\) whenever \(c_k \geq c_*\).
    As \(a\) and \(b\) are coprime, this implies every \(c_k\) between \(c_*\) and \(c_* + b-1\) is attained in the cycle.
    Similarly, modulo \(a\) the value of \(c_k - c_*\) does not change in successful iterations, but does change by \(+b\) whenever \(c_k < c_*\) and so \(c_k\) will attain each value between \(c_* - a\) and \(c_* - 1\).
    Overall, the cycle starting at \(c_0 = c_*\) includes \(c_*-a, \dots, c_*-1\) and \(c_*, \dots, c_*+b-1\) and so it has length \(a+b\).
    As the cycle includes all values in \([c_* - a, c_* + b)\), even if \(c_0 \neq c_*\) there will eventually be an iteration with \(c_k\) in the given range, which will then cycle through all \(a+b\) values.
    During that cycle \(a\) iterations are unsuccessful and \(b\) iterations are successful, meaning that the ratio of unsuccessful to successful iterations is exactly \(\alpha = \frac{a}{b}\) as claimed.
\end{proof}

\section{Numerical Illustration}\label{sec:numerical-illustration}

\begin{figure}
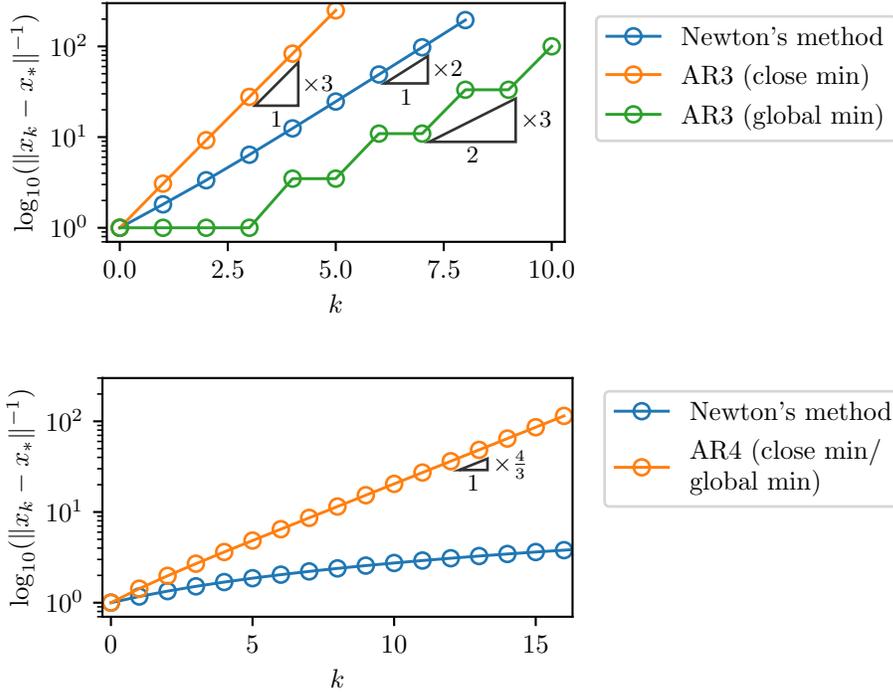
\label{fig:convergence-illustration}
    \centering
    \subimport{./plots}{convergence_pos_def_example.pgf}
    \subimport{./plots}{convergence_sharpness_example.pgf}
    \vspace{-1em}
    \caption{
        Convergence of Newton's method and AR\(p\) for \cref{exm:sigma-oscillations} (top) and \cref{exm:sharp-pq-convergence} (bottom, \(p = 4\), \(q = 4\)).
        The slope of the triangles correspond to the theoretical order of convergence.
        Variants using the global model minimizer and the model minimizer close to the expansion point (see \cref{thm:close-minimizer}) are shown.
    }
\end{figure}

Since the impact of the abstract convergence orders derived in the previous sections may not become immediately clear, we provide illustrations for the example functions from \cref{exm:sigma-oscillations,exm:sharp-pq-convergence}.
The experiments are run using arbitrary precision arithmetic and are terminated when the distance to the true minimizer is less than $10^{-100}$.
The parameter values of \cref{alg:arp} are chosen as \(\eta = \frac{1}{2}\), \(\gamma_1 = \frac{1}{2}\), \(\gamma_2 = 2\), \(\theta = 0\), \(\sigma_0 = \frac{1}{2}\), \(x_0 = x_* + \frac{1}{10}\).
To showcase convergence orders larger than one, the graphs in \cref{fig:convergence-illustration} show \(\log_{10}(\norm{x_k - x_*}^{-1})\), an estimate of the number of correct digits, on a logarithmic scale.

For the non-degenerate case in \cref{exm:sigma-oscillations}, Newton's method converges quadratically.
This means the number of correct digits roughly doubles at each iteration and thus, on a logarithmic scale, the corresponding line has a slope of \(\log(2)\).
The variant of AR\(3\) using the successful model minimizer that is close to the expansion point, achieves cubic convergence, which corresponds to a line with slope \(\log(3)\).
When the global model minimizer is chosen at each step instead, the fourth iteration is the first successful one, since we have \(\sigma_k > 2\) for the first time.
Afterwards, the value of \(\sigma_k\) oscillates, and we can see that in effect only every second iteration is successful as predicted.
With cubic progress during successful steps and no progress during unsuccessful steps, the overall slope is \(\log(3)/2 = \log(\sqrt{3})\).

The second graph covers \cref{exm:sharp-pq-convergence} with \(p = 4\) and \(q = 4\).
On this degenerate example, Newton's method only achieves linear convergence, which corresponds to a line with decreasing slope.
The fourth-order method, on the other hand, achieves superlinear convergence of order \(\frac{p}{q-1}\), i.e., with slope \(\log(\frac{4}{3})\).

\section{Discussion and Conclusion}\label{sec:discussion}

To round up the results in this paper, we compare them to those of Doikov and Nesterov~\cite{doikov_local_2022} and Cartis, Gould and Toint~\cite[Section 5.3]{cartis_evaluation_2022} and discuss connections to a characterization of local minimizers introduced by Cartis et al.~\cite{cartis_efficient_2025}.

As pointed out in the introduction, our paper can be seen as an extension to the analysis in \cite{doikov_local_2022} for objective functions without a non-smooth component.
We replace a constant regularization parameter that requires knowledge of the Lipschitz constant with an adaptive one, relax global uniform convexity to local uniform convexity, and avoid the need for the exact global model minimizer by requiring only an approximate local one.
This relaxation of assumptions comes at the cost of slightly weaker results.
Compare \cref{thm:function-convergence-pq-rate,thm:gradient-convergence-pq-rate} showing that the function value gap \(f(\vek{x}_k) - f(\vek{x}_*)\) and gradient norm \(\norm{\nabla f(\vek{x}_k)}\) decrease at a \(\frac{p}{q-1}\)th-order rate during successful iterations to Theorems 1 and 2 in \cite{doikov_local_2022}.
The conclusions are the same up to small differences in the constants and the fact that \cref{thm:function-convergence-pq-rate} only applies when the norm of the gradient is small enough, which Doikov and Nesterov can avoid by exploiting convexity of the local model.
Similarly, they do not need to consider the possibility of multiple local minimizers or unsuccessful iterations as in \cref{exm:sigma-oscillations}.

Applying our results to their case, we take \(\sigma_0 \geq \frac{p L_p}{(p+1)!}\), \(\gamma_1 = 1\), \(\theta = 0\) and \(\eta = \frac{1}{2}\).
By \cref{thm:sigma-max} every iteration is successful and so the regularization parameter stays constant and large enough.
Moreover, the convexity of the model implies the sublevel set \(\mathcal{L}_k\) from \cref{thm:connected-component-minimizer} is convex with only one connected component, implying that the global minimizer of the model lies in it and thus by \cref{thm:convergence-pq-rate-good-min} the same asymptotic Q-superlinear convergence of function values and gradients as in \cite{doikov_local_2022} is proved.

Section 5.3 in the book by Cartis, Gould and Toint \cite{cartis_evaluation_2022} considers a wider class of functions, namely ones that are measure-dominated.
The relevant class to compare with is the class of gradient-dominated functions of level \(\mu = \frac{q}{q-1}\), which are a generalization of uniformly convex functions of order \(q\) that satisfy only \cref{eqn:uniform-convexity-function-bounded-by-gradient}, i.e.
\[ f(\vek{x}) - f(\vek{x}_*) \leq \kappa \norm{\nabla f(\vek{x})}^\mu \]
for some \(\kappa > 0\).
This inequality is also sometimes referred to as the Polyak--Łojasiewicz property.
The algorithm analysed in \cite{cartis_evaluation_2022}\footnote{This is Algorithm 4.1.1 in \cite{cartis_evaluation_2022} with \(\e_2 = \e_3 = \infty\).} requires an approximate minimizer of the model in the sense of \cref{eqn:approximate-minimizer} at each iteration and updates the regularization parameter \(\sigma_k\) in a very similar way to \cref{alg:arp}.
There is some additional flexibility in \cite{cartis_evaluation_2022} by introducing unsuccessful, successful and very successful iterations depending on the value of \(\rho_k\) and by allowing the increase and decrease factors to vary between iterations within given bounds.
Simplifying these updates as done in \cref{alg:arp} does not fundamentally change the behaviour of the method, but was adopted here so that \cref{thm:convergence-pq-rate,thm:oscillations} are easier to state.
A key change, however, is that the algorithm in \cite{cartis_evaluation_2022} enforces a minimum regularization value \(\sigma_{\min} > 0\).
By using this restriction, it is possible to show a lower bound on the function decrease in terms of the gradient norm during successful iterations with an exponent of \(\nu = \frac{p+1}{p}\):
\[
    f(\vek{x}_k) - f(\vek{x}_{k+1}) \geq \eta \sigma_{\min} R_{\max}^{-\frac{p+1}{p}} \norm{\nabla f(\vek{x}_{k+1})}^{\frac{p+1}{p}}
\]

Combining the gradient-domination of level \(\mu\) and the function decrease result with exponent \(\nu\), Cartis, Gould and Toint are able to derive sublinear rates when \(\nu > \mu\), linear rates when \(\nu = \mu\) and superlinear rates for \(\nu < \mu\).
Specifically, for the \(p > q-1\) case considered in this paper we have \(\nu < \mu\), and so they show that only \(O(\log(\log (\e^{-1})))\) iterations are required to achieve a tolerance of \(\e > 0\) for the function value gap \(f(\vek{x}_k) - f(\vek{x}_*)\).
The order of convergence for the function values is not stated explicitly in \cite[Theorem 5.3.3]{cartis_evaluation_2022}, however from the proof we can see that a scaled function value gap \(\omega_k\) satisfies \(\omega_{k+1} \leq \omega_k^{\mu / \nu}\) for all \(k\) large enough if iteration \(k\) is successful \cite[p. 465]{cartis_evaluation_2022}.
This corresponds to Q-convergence of order
\[
    \frac{\mu}{\nu} = \frac{pq}{(p+1)(q-1)} < \frac{p}{q-1}
\]
for the function values during successful iterations.
The proof shows R-convergence of the same order for the gradient norms during successful iterations as well.
The overall order of convergence is again not considered, because it does not change the log-log dependency on the tolerance \(\e\).
With the same setup as in \cref{thm:convergence-pq-rate} one would get R-convergence of order \(\sqrt[1 + \alpha]{\frac{pq}{(p+1)(q-1)}}\) for function values and gradient norms.

With their weaker assumptions on the function class and the addition of \(\sigma_{\min}\), the authors are able to derive impressive results for the \(p < q-1\) and \(p = q-1\) cases, which go beyond the scope of this paper.
For uniformly convex functions of order \(q\) and \(p > q-1\), however, their rates are suboptimal and as shown in \cref{thm:convergence-pq-rate-successful} the R-convergence of gradient norms can be strengthened to Q-convergence.
Note that adding \(\sigma_{\min}\) to \cref{alg:arp} does not change any of the results in \cref{sec:local-convergence}, and only requires \(\sigma_{\min} < \gamma_1 \sigma_*\) in \cref{thm:oscillations}.

There are also connections of the notion of ``right'' and ``wrong'' minimizers in this paper to persistent and transient model minimizers introduced by Cartis et al.~\cite{cartis_efficient_2025}.
The authors of that study observe that the minimizers of regularized third-order Taylor expansions in general do not vary continuously with the regularization parameter \(\sigma\), unlike the global minimizer of regularized first- and second-order polynomials.
A local minimizer \(\vek{y}_k\) of \(m_k\) is persistent if there exists a continuous curve \(\vek{\psi} \colon [\sigma_k , \infty) \to \R^n\) with \(\vek{\psi}(\sigma_k) = \vek{y}_k\) mapping each regularization parameter value to a local minimizer of the correspondingly regularized model, i.e., if one can continuously trace the local minimizer as \(\sigma\) grows to infinity while the Taylor polynomial stays unchanged.
Any local minimizer that is not persistent is called transient.
Cartis et al.~\cite{cartis_efficient_2025} prove that any global minimizer of the AR1 and AR2 subproblem is persistent and that this property breaks for \(p \geq 3\).
Moreover, they develop an approximate way of detecting persistence and observe empirically that the performance of third-order methods can be improved by rejecting points that are deemed approximately transient.
While the exact relationship between persistent and asymptotically successful minimizers remains an open question, we conjecture that for \(\vek{x}_k\) close enough to \(\vek{x}_*\) there always exists a persistent minimizer and that any persistent minimizer lies in the same connected component of the sublevel set of \(m_k\) as \(\vek{x}_k\) and thus leads to a successful iteration according to \cref{thm:connected-component-minimizer}.
If correct, the persistent minimizer definition distinguishes useful from spurious model minimizers both in the local and non-local regime.
This would encourage research into subproblem solvers aiming for (approximations of) persistent instead of global model minimizers.

Overall, in this paper we show that the results of Doikov and Nesterov \cite{doikov_local_2022} for convex functions hold more generally and that those of Cartis, Gould and Toint \cite{cartis_evaluation_2022} can be strengthened for locally uniformly convex functions to \(p / (q-1)\)th-order convergence during successful iterations (\cref{thm:convergence-pq-rate-successful}), but not further (\cref{exm:sharp-pq-convergence}).
Moreover, we clarify the importance of choosing the right local model minimizer for the asymptotic rates as shown in \cref{exm:sigma-oscillations} as well as \cref{thm:convergence-pq-rate,thm:convergence-pq-rate-good-min}.
Two different characterizations of ``right'' are provided (\cref{thm:close-minimizer,thm:connected-component-minimizer}) and an analysis that shows that oscillations of \(\sigma_k\) around some \(\sigma_*\) are fundamentally unavoidable for fully adaptive regularization methods and odd \(p\) if the ``wrong'' minimizer is chosen (\cref{thm:oscillations}).
We hereby elucidate both the advantages and potential pitfalls of higher-order methods from the perspective of local convergence.



\bibliographystyle{siamplain}
\bibliography{citations}

\end{document}